# ON THE RELATION OF DELAY EQUATIONS TO FIRST-ORDER HYPERBOLIC PARTIAL DIFFERENTIAL EQUATIONS

Iasson Karafyllis[*] and Miroslav Krstic[**]

[*]Department of Mathematics, National Technical University of Athens,
Zografou Campus, 15780, Athens, Greece (email: iasonkar@central.ntua.gr )

[**]Dept. of Mechanical and Aerospace Eng., University of California,
San Diego, La Jolla, CA 92093-0411, U.S.A. (email: krstic@ucsd.edu )

**Abstract**
This paper establishes the equivalence between systems described by a single first-order hyperbolic partial differential equation and systems described by integral delay equations. System-theoretic results are provided for both classes of systems (among them converse Lyapunov results). The proposed framework can allow the study of discontinuous solutions for nonlinear systems described by a single first-order hyperbolic partial differential equation under the effect of measurable inputs acting on the boundary and/or on the differential equation. An illustrative example shows that the conversion of a system described by a single first-order hyperbolic partial differential equation to an integral delay system can simplify considerably the solution of the corresponding robust feedback stabilization problem.

## 1. Introduction

The relation of first-order hyperbolic Partial Differential Equations (PDEs) with delay equations is well known. In many cases, a system of First-Order Hyperbolic PDEs (FOH-PDEs) can be transformed to a system described by Retarded Functional Differential Equations (see [8,9,23]), On the other hand, the recent works [12,13,14] have shown that systems of first-order hyperbolic PDEs can be utilized for the stabilization of delay systems.

However, recent works on the control of systems described by FOH-PDEs (see [1,2,3,4,5,12,19,21,22,26,27,28]) have led to the study of non-standard systems of FOH-PDEs with the following features: (i) the boundary conditions are given by functionals of the whole state profile, and (ii) the differential equations involve functionals of the whole state profile and not only point values of the states. It should be noted that for such systems there is no available existence-uniqueness theory similar to the theory of standard FOH-PDEs (see [15,24]): the researchers have utilized transformation arguments which can guarantee important system-theoretic properties. Moreover, the existence of functionals of the whole state profile in the mathematical description of such systems do not allow the answer to the following question:

"What class of delay systems can be used for the description of
non-standard systems of FOH-PDEs?"

In this work, we answer the above question for a class of systems described by a single FOH-PDE. First, we show that the appropriate class of delay systems is a class which has been rarely studied, i.e., systems of the form:

$$x(t) = f(x_t, w_t)$$
$$x(t) \in \Re^n, w(t) \in W \subseteq \Re^m \quad (1.1)$$

where $W \subseteq \Re^m$ is a non-empty locally compact set with $0 \in U$, $r > 0$, $x_t \in L^\infty([-r,0);\Re^n)$, $w_t \in L^\infty([-r,0];W)$ are defined by $(x_t)(s) = x(t+s)$ for $s \in [-r,0)$, $(w_t)(s) = w(t+s)$ for $s \in [-r,0]$ and $f : L^\infty([-r,0);\Re^n) \times L^\infty([-r,0];W) \to \Re^n$ is a mapping. We call the above class of delay systems a system described by Integral Delay Equations (IDEs). Systems of the form (1.1) have been studied in [9,16,17,18,20]. It should be noted at this point that systems described by IDEs have been utilized for a long time in the stabilization of finite-dimensional control systems with input delays (see [11,13,14]): every predictor feedback is a system described by an IDE, whose input is the state of the finite-dimensional control system. Indeed, a predictor feedback for the finite-dimensional system

$$\dot{y}(t) = g(y(t), x(t-r) + v(t))$$
$$y(t) \in \Re^k, x(t) \in \Re^n, v(t) \in \Re^n$$

where we use $x$ to denote the control input and $v \in L^\infty_{loc}(\Re_+;\Re^n)$ denotes the control actuator error, results in a static feedback law of the form:

$$x(t) = f(x_t, y(t) + u(t))$$

where $u \in L^\infty_{loc}(\Re_+;\Re^n)$ denotes the measurement error. Clearly, the above closed-loop system is the feedback interconnection of an time-invariant system described by ODEs and an integral delay system of the form (1.1). Therefore, the study of systems described by IDEs is important on its own. Important system-theoretic properties for the class of systems described by IDEs are provided in Section 2 of the present paper.

Secondly, we show that a class of systems described by a non-standard single FOH-PDE, namely systems of the form:

$$\frac{\partial x}{\partial t}(t,z) + c\frac{\partial x}{\partial z}(t,z) = a(p(t),z)x(t,z) + g(z)p(t), \text{ for } t > 0, z \in (0,1) \quad (1.2)$$

$$p_i(t) = K_i(w(t), x_t), \text{ for } t > 0, i = 1,...,N \text{ and } p(t) = (p_1(t),...,p_N(t))' \in \Re^N \quad (1.3)$$

$$x(t,0) = G(w(t), x_t), \text{ for } t \geq 0 \quad (1.4)$$

where $N > 0$ is an integer, $x(t,z) \in \Re$, $w(t) \in W \subseteq \Re^m$, $W \subseteq \Re^m$ is a locally compact set, $c > 0$ is a constant, $g(z) = (g_1(z),...,g_N(z))$, $G : W \times L^\infty((0,1];\Re) \to \Re$, $K_i : W \times L^\infty((0,1];\Re) \to \Re$ ($i = 1,...,N$) are functionals, $a(p,z)$, $g_i(z)$ ($i = 1,...,N$) are sufficiently regular scalar functions and $x_t$ is the state profile (i.e., $(x_t)(z) = x(t,z)$ for all $t \geq 0$ and $z \in [0,1]$), are equivalent to systems described by IDEs.

The equivalence allows the development of important system-theoretic results for the class of systems described by (1.2), (1.3) and (1.4) (Section 3 of the present paper). The obtained results can allow the study of measurable inputs $w(t) \in W \subseteq \Re^m$, which is an important feature because in many cases boundary conditions of the form (1.4) result from the implementation of stabilizing feedback laws (see [1,2,3,4,5,12,21,22]). However, the control action comes together with control actuator errors and measurement errors, which are typically modeled by measurable inputs.

Finally, the obtained results can allow the development of a methodology for the solution of control problems for systems described by a non-standard single FOH-PDE, namely systems of the form (1.2), (1.3), (1.4). The methodology is presented in Section 4 of the present paper by means of an illustrative example. The example shows many important features of the proposed methodology and a comparison is made with other existing methodologies.

**Notation.** Throughout this paper, we adopt the following notation:

* $\Re_+ := [0, +\infty)$.
* Let $I \subseteq \Re$ be an interval and $U \subseteq \Re^m$ be a set. By $L^\infty(I;U)$ we denote the space of measurable and essentially bounded functions $u(\cdot)$ defined on $I$ and taking values in $U \subseteq \Re^m$. By $L^\mu(I;U)$, where $\mu \in [1,+\infty)$, we denote the space of measurable functions $u(\cdot)$ defined on $I$ and taking values in $U \subseteq \Re^m$ for which $\int_I |u(s)|^\mu ds < +\infty$. By $L^\infty_{loc}(I;U)$ we denote the space of measurable and locally essentially bounded functions $u(\cdot)$ defined on $I$ and taking values in $U \subseteq \Re^m$, i.e., for every compact interval $J \subseteq I$ it holds that $u \in L^\infty(J;U)$. For every $x \in L^\infty(I;\Re^n)$, where $I \subseteq \Re$ is an interval, we denote by $\|x\|$ the essential supremum of $x$, i.e., $\|x\| = \sup_{s \in I} |x(s)|$. For $x \in L^\infty(I;\Re^n)$, where $I \subseteq \Re$ is an interval with $[a,b) \subseteq I$, $a < b$, we define $\|x\|_{[a,b)} = \sup_{a \leq s < b} |x(s)|$.
* For a vector $x \in \Re^n$, we denote by $x'$ its transpose and by $|x|$ its Euclidean norm.

## 2. Systems Described by Integral Delay Equations

We consider system (1.1) under the following assumptions:

**(H1)** *There exist non-decreasing functions $a : \Re_+ \to \Re_+$, $M : \Re_+ \to \Re_+$, $N : \Re_+ \to \Re_+$ such that for every $R > 0$ and for every $x, y \in L^\infty([-r,0);\Re^n)$, $w \in L^\infty([-r,0];W)$ with $\|x\| \leq R$, $\|y\| \leq R$, $\|w\| \leq R$ the following inequalities hold:*

$$|f(x,w) - f(y,w)| \leq N(R) h \sup_{-h \leq s < 0} |x(s) - y(s)| + M(R) \sup_{-r \leq s < -h} |x(s) - y(s)|, \text{ for all } h \in (0,r) \quad (2.1)$$

$$|f(x,w)| \leq a(R) \quad (2.2)$$

**(H2)** *For every $\delta > 0$, $x \in L^\infty([-r,\delta);\Re^n)$, $w \in L^\infty([-r,\delta];W)$, the function $z : [-r,\delta) \to \Re^n$ defined by $z(t) = x(t)$ for $t \in [-r,0)$ a.e. and $z(t) = f(x_t, w_t)$ for $t \in [0,\delta)$ a.e. satisfies $z \in L^\infty([-r,\delta);\Re^n)$.*

Systems of the form (1.1) can represent integral delay equations of the form:

$$x(t) = \sum_{i=1}^{N} p_i(x(t-\tau_i), w(t-\tau_i)) + \int_0^r q(s, x(t-s), w(t-s)) ds$$

where $r > 0$, $0 < \tau_1 \leq \tau_2 \leq \ldots \leq \tau_N \leq r$ are constants, $p_i : \Re^n \times W \to \Re^n$ ($i = 1, \ldots, N$) are locally Lipschitz mappings, $W \subseteq \Re^m$ is a non-empty locally compact set and $q : \Re \times \Re^n \times W \to \Re^n$ is a locally Lipschitz mapping. It can be shown that assumptions (H1), (H2) hold for this class of systems. The fact that assumption (H2) holds for this class of systems is a direct consequence of Lemma 1 on page 4 in [6].

The following result shows us that system (1.1) can be regarded as a time-invariant system with inputs, whose state space is $L^\infty([-r,0];\Re^n)$ and whose input set is $L^\infty([-r,0];W)$. System (1.1) satisfies the classical semigroup property, the Boundedness-Implies-Continuation (BIC) property and the property of Lipschitz dependence on the initial conditions (see [10]).

**Theorem 2.1:** *Consider system (1.1) under assumptions (H1), (H2). Then for every $x_0 \in L^\infty([-r,0];\Re^n)$, $w \in L^\infty_{loc}([-r,+\infty);W)$ there exists $t_{\max} = t_{\max}(x_0, w) \in (0, +\infty]$ with*

$$t_{\max} > \frac{r}{1 + 2rN\left(5a\left(\max\left(\|x_0\|, \sup_{-r \leq s \leq r} |w(s)|\right)\right)\right)}$$

*and a unique mapping $x \in L^\infty_{loc}([-r, t_{\max}); \Re^n)$ satisfying $x(t) = f(x_t, w_t)$ for $t \in [0, t_{\max})$ a.e. and $x(t) = x_0(t)$ for $t \in [-r,0)$ a.e.. Moreover, if $t_{\max} < +\infty$ then $\limsup_{t \to t_{\max}^-} \|x_t\| = +\infty$. Finally, there exist non-decreasing functions $P : \Re_+ \to \Re_+$, $G : \Re_+ \to \Re_+$ such that for every $x_0 \in L^\infty([-r,0];\Re^n)$, $y_0 \in L^\infty([-r,0];\Re^n)$, $w \in L^\infty_{loc}([-r,+\infty);W)$ it holds that:*

$$\|x_t - y_t\| \leq G(s(t))\exp(P(s(t))t)\|x_0 - y_0\|, \text{ for all } t \in [0, \delta) \quad (2.3)$$

*where $\delta = \min(t_{\max}(x_0, w), t_{\max}(y_0, w))$, $s(t) := \max\left(\sup_{0 \leq \tau \leq t}\|x_\tau\|, \sup_{0 \leq \tau \leq t}\|y_\tau\|, \sup_{0 \leq \tau \leq t}\|w_\tau\|\right)$ and $x \in L^\infty_{loc}([-r, t_{\max}(x_0, w));\Re^n)$, $y \in L^\infty_{loc}([-r, t_{\max}(y_0, w));\Re^n)$ are the unique mappings satisfying $x(t) = f(x_t, w_t)$ for $t \in [0, t_{\max}(x_0, w))$ a.e., $y(t) = f(y_t, w_t)$ for $t \in [0, t_{\max}(y_0, u))$ a.e. and $x(t) = x_0(t)$, $y(t) = y_0(t)$ for $t \in [-r, 0)$ a.e..*

**Proof:** Let arbitrary $x_0 \in L^\infty([-r,0];\Re^n)$, $w \in L^\infty_{loc}([-r,+\infty);W)$. Define $s := \max\left(\|x_0\|, \sup_{-r \leq t \leq r}|u(t)|\right)$ and select arbitrary $p \in \Re^n$ with $|p| \leq s$. Define $x^{(1)}(t) = x_0(t)$ for $t \in [-r, 0)$ a.e. and $x^{(1)}(t) = p$ for $t \in [0, r)$. Clearly, $x^{(1)} \in L^\infty([-r, r); \Re^n)$ with $\sup_{0 \leq \tau < r}\|x_\tau^{(1)}\| \leq s$.

Without loss of generality we may assume that the non-decreasing function $a : \Re_+ \to \Re_+$ involved in (2.2) satisfies $a(t) \geq t$ for all $t \geq 0$. Let $\delta \in (0, r)$ be such that:

$$2N(R)\delta < 1, \text{ for } R = 5a(s) \quad (2.4)$$

Using induction and assumption (H2), we can guarantee that the sequence of functions $x^{(k)} : [-r, \delta) \to \Re^n$ for $k \geq 1$ defined by:

$x^{(k)}(t) = x_0(t)$ for $t \in [-r, 0)$ a.e. and $x^{(k)}(t) = f(x_t^{(k-1)}, w_t)$ for $t \in [0, \delta)$ a.e. for $k \geq 2$ satisfies $x^{(k)} \in L^\infty\left([-r, \delta); \Re^n\right)$ for $k \geq 1$.

Notice that (2.2) in conjunction with the facts that $\sup_{0 \leq \tau < r} \|x_\tau^{(1)}\| \leq s$ and $\sup_{0 \leq \tau \leq r} \|u_\tau\| \leq s$ implies that

$$\sup_{0 \leq \tau < \delta} \|x_\tau^{(1)}\| \leq a(s) \text{ and } \sup_{0 \leq \tau < \delta} \|x_\tau^{(2)}\| \leq a(s) \tag{2.5}$$

We next claim that

$$\max_{l=1,\ldots,k}\left(\sup_{0 \leq \tau < \delta} \|x_\tau^{(l)}\|\right) \leq 5a(s) \text{ for all } k \geq 1 \tag{2.6}$$

The proof of the claim is made by induction. We first notice that by virtue of (2.5) the claim is true for $k=1$ and $k=2$. Next suppose that the claim is true for certain $k \geq 2$. Notice that, by virtue of (2.1) and the fact that $R = 5a(s)$, the following inequality holds for all $l = 2, \ldots, k$ and for almost all $t \in [0, \delta)$:

$$\begin{aligned}
\left|x^{(l+1)}(t) - x^{(l)}(t)\right| &= \left|f(x_t^{(l)}, w_t) - f(x_t^{(l-1)}, w_t)\right| \\
&\leq N(R)\delta \sup_{-\delta \leq s < 0}\left|(x_t^{(l)})(s) - (x_t^{(l-1)})(s)\right| + M(R) \sup_{-r \leq s < -\delta}\left|(x_t^{(l)})(s) - (x_t^{(l-1)})(s)\right| \\
&= N(R)\delta \sup_{-\delta \leq s < 0}\left|x^{(l)}(t+s) - x^{(l-1)}(t+s)\right| + M(R) \sup_{-r \leq s < -\delta}\left|x^{(l)}(t+s) - x^{(l-1)}(t+s)\right| \\
&= N(R)\delta \left\|x^{(l)} - x^{(l-1)}\right\|_{[-r, \delta)}
\end{aligned} \tag{2.7}$$

Since (2.7) holds for almost all $t \in [0, \delta)$ and since $x^{(l+1)}(t) = x^{(l)}(t) = x_0(t)$ for almost all $t \in [-r, 0)$ and for all $l \geq 1$, we get:

$$\left\|x^{(l+1)} - x^{(l)}\right\|_{[-r, \delta)} \leq N(R)\delta \left\|x^{(l)} - x^{(l-1)}\right\|_{[-r, \delta)}, \text{ for all } l = 2, \ldots, k \tag{2.8}$$

Inequality (2.8) implies that:

$$\left\|x^{(2+m)} - x^{(1+m)}\right\|_{[-r, \delta)} \leq (N(R)\delta)^m \left\|x^{(2)} - x^{(1)}\right\|_{[-r, \delta)}, \text{ for all } m = 0, \ldots, k-1 \tag{2.9}$$

Using the triangle inequality repeatedly and the fact that $N(R)\delta \leq 1/2$ (a consequence of (2.4)), we get for all $l = 2, \ldots, k$:

$$\begin{aligned}
\left\|x^{(l+1)} - x^{(1)}\right\|_{[-r, \delta)} &\leq \sum_{m=0}^{l-1} \left\|x^{(m+2)} - x^{(m+1)}\right\|_{[-r, \delta)} \\
&\leq \left\|x^{(2)} - x^{(1)}\right\|_{[-r, \delta)} \sum_{m=0}^{l-1} (N(R)\delta)^m = \left\|x^{(2)} - x^{(1)}\right\|_{[-r, \delta)} \frac{1 - (N(R)\delta)^l}{1 - N(R)\delta} \leq 2\left\|x^{(2)} - x^{(1)}\right\|_{[-r, \delta)}
\end{aligned} \tag{2.10}$$

It follows from (2.5) and (2.10) that

$$\sup_{0\leq \tau <\delta}\left\|x_\tau^{(k+1)}\right\| = \left\|x^{(k+1)}\right\|_{[-r,\delta)} \leq 5a(s) \tag{2.11}$$

and therefore we have proved (2.6).

The existence of $x \in L^\infty\big([-r,\delta);\Re^n\big)$ satisfying $x(t) = f(x_t, w_t)$ for $t \in [0,\delta)$ a.e. and $x(t) = x_0(t)$ for $t \in [-r,0)$ a.e. follows closely the proof of Banach's fixed point theorem on the closed set $L^\infty\big([-r,\delta);B\big)$, where $B \subseteq \Re^n$ denotes the closed ball of radius $5a(s)$. Since (2.9) holds for all $m \geq 0$, the sequence $x^{(k)} \in L^\infty\big([-r,\delta);B\big)$ is a Cauchy sequence. Therefore, there exists a unique limit $x \in L^\infty\big([-r,\delta);B\big)$. Assumption (H1) (inequality (2.1)) implies that the right hand side of the inequality

$$\left|x(t) - f(x_t, w_t)\right| \leq \left|x(t) - x^{(k)}(t)\right| + \left|f(x_t, w_t) - f(x_t^{(k-1)}, w_t)\right|$$

can be made arbitrarily small for sufficiently large $k \geq 1$ and for almost all $t \in [0,\delta)$ and consequently, we obtain $x(t) = f(x_t, w_t)$ for $t \in [0,\delta)$ a.e..

In order to prove uniqueness, we use a contradiction argument. Suppose that there exists $y \in L^\infty\big([-r,\delta);\Re^n\big)$ satisfying $y(t) = f(y_t, w_t)$ for $t \in [0,\delta)$ a.e. and $y(t) = x_0(t)$ for $t \in [-r,0)$ a.e. with $\sup_{-r \leq s < \delta}|y(s) - x(s)| > 0$. Let $p \geq 0$ be the least upper bound of all $t \in [0,\delta)$ with $\sup_{-r \leq s < t}|y(s) - x(s)| = 0$. A contradiction argument shows that $p < \delta$. Since $p \in [0,\delta)$ is the least upper bound of all $t \in [0,\delta)$ with $\sup_{-r \leq s < t}|y(s) - x(s)| = 0$, it follows that $\sup_{-r \leq s < t}|y(s) - x(s)| > 0$ for all $t \in (p,\delta)$. A contradiction argument shows that $\sup_{-r \leq s < p}|y(s) - x(s)| = 0$. Define $R = \max\left(\sup_{-r \leq s < \delta}|y(s)|, \sup_{-r \leq s < \delta}|x(s)|, \sup_{-r \leq s \leq \delta}|w(s)|\right)$. It follows from (2.1) for almost all $t \in (p,\delta)$ and $h = t - p$:

$$|x(t) - y(t)| = |f(x_t, w_t) - f(y_t, w_t)| \leq N(R)(t-p)\sup_{p \leq s < t}|x(s) - y(s)| \tag{2.12}$$

Inequality (2.12) implies $\sup_{p < s < t}|x(s) - y(s)| \leq N(R)(t-p)\sup_{p \leq s < t}|x(s) - y(s)|$ for all $t \in (p,\delta)$ and since $\sup_{p < s < t}|x(s) - y(s)| = \sup_{p \leq s < t}|x(s) - y(s)|$, we obtain $\sup_{p \leq s < t}|x(s) - y(s)| = 0$ for all $t < \min\left(\delta, p + \frac{1}{N(R)+1}\right)$, i.e., $\sup_{-r \leq s < t}|y(s) - x(s)| = 0$ for all $t < \min\left(\delta, p + \frac{1}{N(R)+1}\right)$, a contradiction.

Since the unique function $x \in L^\infty\big([-r,\delta);\Re^n\big)$ satisfying $x(t) = f(x_t, w_t)$ for $t \in [0,\delta)$ a.e. and $x(t) = x_0(t)$ for $t \in [-r,0)$ a.e. is bounded, it follows (by repeating the arguments above) that there exists $\delta' > \delta$ and a unique function $x \in L^\infty\big([-r,\delta');\Re^n\big)$ satisfying $x(t) = f(x_t, w_t)$ for $t \in [0,\delta')$ a.e. and $x(t) = x_0(t)$ for $t \in [-r,0)$ a.e.. More specifically, by virtue of (2.4), $\delta' > \delta$ satisfies

$$\delta + r \geq \delta' \geq \delta + \frac{r}{1 + 2rN\left(5a\left(\max\left(\sup_{-r \leq s < \delta}|x(s)|, \sup_{-r \leq s \leq \delta+r}|w(s)|\right)\right)\right)} \tag{2.13}$$

Let $t_{max} = t_{max}(x_0, w) \in (0, +\infty]$ be the least upper bound of all $\delta > 0$ for which there exists function $x \in L^\infty([-r, \delta); \Re^n)$ satisfying $x(t) = f(x_t, w_t)$ for $t \in [0, \delta)$ a.e. and $x(t) = x_0(t)$ for $t \in [-r, 0)$ a.e.. If $t_{max} < +\infty$ (i.e., is finite) then for every $\varepsilon > 0$ there exists $\delta > t_{max} - \varepsilon$ and a function $x \in L^\infty([-r, \delta); \Re^n)$ satisfying $x(t) = f(x_t, w_t)$ for $t \in [0, \delta)$ a.e. and $x(t) = x_0(t)$ for $t \in [-r, 0)$ a.e.. It follows in any case ($t_{max} < +\infty$ or $t_{max} = +\infty$) that there exists a function $x \in L^\infty_{loc}([-r, t_{max}); \Re^n)$ satisfying $x(t) = f(x_t, w_t)$ for $t \in [0, t_{max})$ a.e. and $x(t) = x_0(t)$ for $t \in [-r, 0)$ a.e..

For the case $t_{max} < +\infty$ we notice that (2.13) implies the inequality 
$$t_{max} \geq \delta + \frac{r}{1 + 2rN\left(5a\left(\max\left(\sup_{-r \leq s < \delta} |x(s)|, \sup_{-r \leq s \leq \delta + r} |w(s)|\right)\right)\right)} \text{ or } N\left(5a\left(\max\left(\sup_{-r \leq s < \delta} |x(s)|, \sup_{-r \leq s \leq \delta + r} |w(s)|\right)\right)\right) \geq \frac{r - \varepsilon}{2r\varepsilon}$$
for all $\varepsilon > 0$ and $\delta > t_{max} - \varepsilon$. The inequality $N\left(5a\left(\max\left(\sup_{-r \leq s < \delta} |x(s)|, \sup_{-r \leq s \leq \delta + r} |w(s)|\right)\right)\right) \geq \frac{r - \varepsilon}{2r\varepsilon}$ cannot hold for all $\varepsilon > 0$ if $\sup_{-r \leq s < t_{max}} |x(s)| < +\infty$. Therefore, if $t_{max} < +\infty$ then $\limsup_{t \to t_{max}^-} \|x_t\| = +\infty$.

Finally, we proceed to the constructive proof of estimate (2.3). Let $x_0 \in L^\infty([-r, 0); \Re^n)$, $y_0 \in L^\infty([-r, 0); \Re^n)$, $w \in L^\infty_{loc}([-r, +\infty); W)$ and let $x \in L^\infty_{loc}([-r, t_{max}(x_0, w)); \Re^n)$, $y \in L^\infty_{loc}([-r, t_{max}(y_0, w)); \Re^n)$ be the unique mappings satisfying $x(t) = f(x_t, w_t)$ for $t \in [0, t_{max}(x_0, w))$ a.e., $y(t) = f(y_t, w_t)$ for $t \in [0, t_{max}(y_0, w))$ a.e. and $x(t) = x_0(t)$, $y(t) = y_0(t)$ for $t \in [-r, 0)$ a.e.. Define $\delta = \min(t_{max}(x_0, w), t_{max}(y_0, w))$. For all $t \in (0, \delta)$ we define $s(t) := \max\left(\sup_{0 \leq \tau \leq t} \|x_\tau\|, \sup_{0 \leq \tau \leq t} \|y_\tau\|, \sup_{0 \leq \tau \leq t} \|w_\tau\|\right)$ and $h = \frac{t}{N}$, where $N = 1 + \left[\frac{t(1 + 2rN(s(t)))}{\min(r, t)}\right]$ and $\left[\frac{t(1 + 2rN(s(t)))}{\min(r, t)}\right]$ denotes the integer part of the real $\frac{t(1 + 2rN(s(t)))}{\min(r, t)}$. Inequality (2.1) implies for almost all $v \in [0, h]$:

$$|x(v) - y(v)| = |f(x_v, w_v) - f(y_v, w_v)| \leq N(s(t)) h \sup_{-h \leq q < 0} |x(v + q) - y(v + q)| + M(s(t)) \|x_0 - y_0\| \quad (2.14)$$

Notice that since $N(s(t)) h \leq 1/2$ we get from (2.14):

$$\sup_{-r \leq q < h} |x(q) - y(q)| \leq (1 + 2M(s(t))) \|x_0 - y_0\| \quad (2.15)$$

Using induction and similar arguments we can prove that the following inequality holds for $i = 1, ..., N$

$$\sup_{-r \leq q < ih} |x(q) - y(q)| \leq (1 + 2M(s(t)))^i \|x_0 - y_0\| \quad (2.16)$$

Since $N = 1 + \left[\frac{t(1 + 2rN(s(t)))}{\min(r, t)}\right] \leq 1 + \frac{t(1 + 2rN(s(t)))}{\min(r, t)} \leq 2 + 2rN(s(t)) + \frac{t(1 + 2rN(s(t)))}{r}$, we obtain from (2.16) for $i = N$:

$$\sup_{-r \leq q < t} |x(q) - y(q)| \leq (1 + 2M(s(t)))^{2 + 2rN(s(t))} \exp\left((2N(s(t)) + r^{-1}) \ln(1 + 2M(s(t))) t\right) \|x_0 - y_0\| \quad (2.17)$$

Inequality (2.17) directly implies estimate (2.3) with $G(s) := (1 + 2M(s))^{2+2rN(s)}$ and $P(s) := (2N(s) + r^{-1}) \ln(1 + 2M(s))$.

The proof is complete. ◁

Having proved Theorem 2.1, we are justified for arbitrary $x_0 \in L^\infty([-r,0); \Re^n)$, $w \in L^\infty_{loc}([-r,+\infty); W)$ to name the unique mapping $x \in L^\infty_{loc}([-r,t_{max}); \Re^n)$ satisfying $x(t) = f(x_t, w_t)$ for $t \in [0, t_{max})$ a.e., $x(t) = x_0(t)$ for $t \in [-r, 0)$ a.e. for certain $t_{max} = t_{max}(x_0, w) \in (0, +\infty]$ as "*the solution $x(t)$ of (1.1) with initial condition $x(t) = x_0(t)$ for $t \in [-r, 0)$ a.e., corresponding to input $w \in L^\infty_{loc}([-r,+\infty); W)$*", without explicit reference to the maximal existence time of the solution $t_{max} = t_{max}(x_0, w) \in (0, +\infty]$.

A number of system-theoretic properties can be proved by using the results of Theorem 2.1. More specifically, we study system (1.1) under assumptions (H1), (H2) and the following assumption:

**(H3)** $w = (d, u)$, where $d \in D \subseteq \Re^{m_1}$ *is a compact set*, $u \in U \subseteq \Re^{m_2}$ *is a locally compact set with* $0 \in U$ *and* $m = m_1 + m_2$. *Moreover, there exists* $b \in K_\infty$ *such that for every* $x \in L^\infty([-r,0); \Re^n)$, $d \in L^\infty([-r,0]; D)$, $u \in L^\infty([-r,0]; U)$, *the following inequality holds:*

$$|f(x,d,u)| \leq b\left(\max\left(\|x\|, \sup_{-r \leq s \leq 0} |u(s)|\right)\right) \quad (2.18)$$

Assumption (H3) means that the system's right-hand side is uniformly bounded with respect to the disturbance $d$ and the right-hand side is zero when the state $x$ and input $u$ are zero, irrespective of the disturbance $d$. This means, in particular, that the origin is an equilibrium when the input $u$ is zero, irrespective of the value of the disturbance $d$, which excludes, for example, systems that have an additive disturbance in the model's right-hand side. In Example 2.7 we provide a simple system that captures some of the essence of assumption (H3) and the system's dependence on $d$.

We are ready to give a list of properties for system (1.1) that are derived from the results of Theorem 2.1.

Property 1: Robustness of the equilibrium point.

Assumption (H3) guarantees that $0 \in L^\infty([-r,0); \Re^n)$ is an equilibrium point for system (1.1), when $u$ is zero and for any $d$. However, in order to study the robust stability properties of the equilibrium point of (1.1) we need a stronger assertion, namely that $0 \in L^\infty([-r,0); \Re^n)$ is a robust equilibrium point for the system (1.1) with input $u \in U \subseteq \Re^{m_2}$ (see [10]).

**Theorem 2.2:** *Consider system (1.1) under assumptions (H1), (H2), (H3). The function $0 \in L^\infty([-r,0); \Re^n)$ is a robust equilibrium point for the system (1.1) with input $u \in U \subseteq \Re^{m_2}$, i.e., for all $\varepsilon > 0$ and $T \geq 0$ there exists $\delta := \delta(\varepsilon, T) > 0$ such that for every $x_0 \in L^\infty([-r,0); \Re^n)$, $d \in L^\infty([-r,+\infty); D)$, $u \in L^\infty_{loc}([-r,+\infty); U)$ with $\|x_0\| + \sup_{t \geq 0} \|u_t\| < \delta$ there exists $t_{max} = t_{max}(x_0, d, u) \in (T, +\infty]$ and a unique mapping $x \in L^\infty_{loc}([-r, t_{max}); \Re^n)$ satisfying $x(t) = f(x_t, d_t, u_t)$ for $t \in [0, t_{max})$ a.e., $x(t) = x_0(t)$ for $t \in [-r, 0)$ a.e. and $\|x_t\| < \varepsilon$ for all $t \in [0, T]$.*

**Proof:** Since $D \subseteq \Re^{m_1}$ is compact, we are in the position to define $Q := \max_{d \in D} |d|$. Without loss of generality, we may assume that the function $b \in K_\infty$ involved in (2.18) satisfies $b(s) \geq s$ for all $s \geq 0$. Let arbitrary $\varepsilon > 0$, $T \geq 0$ be given and define:

$$\delta(\varepsilon, T) := \kappa^{-1}(\varepsilon/2), \quad \kappa := g^{(l)} = \underbrace{g \circ ... \circ g}_{l \text{ times}}, \quad l = \left[\frac{T}{\rho}\right] + 1 \text{ and } \rho = \frac{r}{1 + 2rN(5b(Q+\varepsilon))} \quad (2.19)$$

where $g(s) := 5b(s)$ for all $s \geq 0$, $\left[\frac{T}{\rho}\right]$ is the integer part of $\frac{T}{\rho}$ and $N : \Re_+ \to \Re_+$ is the function involved in (2.1). Let arbitrary $x_0 \in L^\infty([-r,0);\Re^n)$, $d \in L^\infty([-r,+\infty);D)$, $u \in L^\infty_{loc}([-r,+\infty);U)$ with $\|x_0\| + \sup_{t \geq 0}\|u_t\| < \delta$ and consider the unique mapping $x \in L^\infty_{loc}([-r,t_{\max});\Re^n)$ satisfying $x(t) = f(x_t, d_t, u_t)$ for $t \in [0, t_{\max})$ a.e., $x(t) = x_0(t)$ for $t \in [-r,0)$ a.e. for certain $t_{\max} = t_{\max}(x_0, d, u) \in (0, +\infty]$. Using (2.18), the fact that $\sup_{t \geq 0}\|u_t\| < \delta$, the fact that $\delta \leq \varepsilon$ (a direct consequence of definition (2.19)), the fact that $Q := \max_{d \in D} |d|$ and proceeding exactly as in the proof of Theorem 2.1, we can show that:

"For all $t_0 \in [0, t_{\max})$ with $\|x_{t_0}\| \leq \varepsilon$ it holds that $\|x_t\| \leq g(\max(\|x_{t_0}\|, \delta))$, for all $t \in [t_0, t_0 + \rho]$" (2.20)

We next claim that $\|x_{i\rho}\| \leq g^{(i)}(\delta)$, for all $i = 1,...,l$ ($\|x_{i\rho}\|$ is the norm $\|x_t\|$ for $t = i\rho$). The claim is a direct consequence of property (2.20) and definition (2.19). Therefore, we get $\|x_t\| \leq \kappa(\delta)$ for all $t \in [0, l\rho]$. Since $l\rho \geq T$ (a direct consequence of the fact that $l = \left[\frac{T}{\rho}\right] + 1$; recall (2.19)) and $\delta = \kappa^{-1}(\varepsilon/2)$ (recall (2.19)), we obtain $\|x_t\| < \varepsilon$ for all $t \in [0, T]$. The proof is complete. ◁

Property 2: Robust Global Asymptotic Stability and Input-to-State Stability

Using the results contained in [10], we are in a position to define the notion of (Uniform) Robust Global Asymptotic Stability and Input-to-State Stability (ISS) for system (1.1) under assumptions (H1), (H2), (H3). The notion of ISS for systems described by IDEs is completely analogous to the corresponding notion introduced by E. D. Sontag for finite-dimensional systems in [25].

**Definition 2.3:** *Consider system (1.1) under assumptions (H1), (H2), (H3) and assume that for every $x_0 \in L^\infty([-r,0);\Re^n)$, $d \in L^\infty([-r,+\infty);D)$, $u \in L^\infty_{loc}([-r,+\infty);U)$, there exists $x \in L^\infty_{loc}([-r,+\infty);\Re^n)$ satisfying $x(t) = f(x_t, d_t, u_t)$ for $t \in [0,+\infty)$ a.e., $x(t) = x_0(t)$ for $t \in [-r,0)$ a.e.. We say that system (1.1) is Input-to-State Stable (ISS) from the input $u \in U \subseteq \Re^{m_2}$ uniformly in $d \in D \subseteq \Re^{m_1}$, if there exists a continuous non-decreasing function $\gamma : \Re_+ \to \Re_+$ such that the following properties hold:*

**Robust Lagrange Stability:** *For every $\varepsilon > 0$ it holds that*

$$\sup\left\{\|x_t\| - \sup_{0 \leq s \leq t} \gamma(\|u_s\|) : t \geq 0, \|x_0\| \leq \varepsilon, d \in L^\infty([-r,+\infty);D), u \in L^\infty_{loc}([-r,+\infty);U)\right\} < +\infty.$$

**Robust Lyapunov Stability:** *For every $\varepsilon > 0$ there exists $\delta := \delta(\varepsilon) > 0$ such that*

$$\sup\left\{\|x_t\| - \sup_{0 \leq s \leq t} \gamma(\|u_s\|) : t \geq 0, \|x_0\| \leq \delta, d \in L^\infty([-r,+\infty);D), u \in L^\infty_{loc}([-r,+\infty);U)\right\} < \varepsilon.$$

**Uniform Robust Attractivity:** *For every $\varepsilon > 0$ and $R \geq 0$ there exists $\tau := \tau(\varepsilon, R) > 0$ such that*

$$\sup\left\{\|x_t\| - \sup_{0 \leq s \leq t} \gamma(\|u_s\|) : t \geq \tau, \|x_0\| \leq R, d \in L^\infty([-r,+\infty);D), u \in L^\infty_{loc}([-r,+\infty);U)\right\} < \varepsilon.$$

*If $U = \{0\}$ then we say that $0 \in L^\infty([-r,0);\Re^n)$ is (Uniformly) Robustly Globally Asymptotically Stable (RGAS) for (1.1).*

Lemma 2.1, Lemma 4.2 in [10] and Theorem 2.2 in [10] give us the following result.

**Theorem 2.4:** *Consider system (1.1) under assumptions (H1), (H2), (H3). Then the following statements are equivalent:*

**(a)** *System (1.1) is ISS from the input $u \in U \subseteq \Re^{m_2}$ uniformly in $d \in D \subseteq \Re^{m_1}$.*

**(b)** *There exists $\sigma \in KL$ and a continuous, non-decreasing function $\gamma : \Re_+ \to \Re_+$ such that for every $x_0 \in L^\infty([-r,0);\Re^n)$, $d \in L^\infty([-r,+\infty);D)$, $u \in L^\infty_{loc}([-r,+\infty);U)$, there exists $x \in L^\infty_{loc}([-r,+\infty);\Re^n)$ satisfying $x(t) = f(x_t, d_t, u_t)$ for $t \in [0,+\infty)$ a.e., $x(t) = x_0(t)$ for $t \in [-r,0)$ a.e. and $\|x_t\| \leq \sigma(\|x_0\|, t) + \sup_{0 \leq s \leq t} \gamma(\|u_s\|)$ for all $t \geq 0$.*

*Furthermore, if $U = \{0\}$ then the following statement is equivalent to statements (a), (b):*

**(c)** *System (1.1) is Robustly Forward Complete, i.e., for every $T \geq 0$ and $R \geq 0$ it holds that $\sup\{\|x_t\| : t \in [0,T], \|x_0\| \leq R, d \in L^\infty([-r,+\infty);D)\} < +\infty$, where $x(t)$ denotes the solution of (1.1) with $x(t) = x_0(t)$ for $t \in [-r,0)$ a.e., corresponding to input $d \in L^\infty([-r,+\infty);D)$, and the property of Uniform Robust Attractivity of Definition 2.3 holds.*

Property 3: Lyapunov Characterization of RGAS

Theorem 3.4 in [10] and the results of Theorems 2.1, 2.2 allow us to obtain a complete Lyapunov characterization for the RGAS property for system (1.1).

**Theorem 2.5:** *Consider system (1.1) under assumptions (H1), (H2), (H3) and assume that $U = \{0\}$. The equilibrium point $0 \in L^\infty([-r,0);\Re^n)$ is RGAS for (1.1) if and only if there exists a functional $V : L^\infty([-r,0);\Re^n) \to \Re_+$, a non-decreasing function $Q : \Re_+ \to \Re_+$ and functions $a_1, a_2 \in K_\infty$ such that the following inequalities hold:*

$$a_1(\|x\|) \leq V(x) \leq a_2(\|x\|), \text{ for all } x \in L^\infty([-r,0);\Re^n) \quad (2.21)$$

$$|V(x) - V(y)| \leq Q(\max(\|x\|, \|y\|))\|x - y\|, \text{ for all } x, y \in L^\infty([-r,0);\Re^n) \quad (2.22)$$

*Moreover, for every $x_0 \in L^\infty([-r,0);\Re^n)$, $d \in L^\infty([-r,+\infty);D)$, the solution $x(t)$ of (1.1) with initial condition $x(t) = x_0(t)$ for $t \in [-r,0)$ a.e., corresponding to input $d \in L^\infty([-r,+\infty);D)$ satisfies:*

$$V(x_t) \leq \exp(-t)V(x_0), \text{ for all } t \geq 0 \quad (2.23)$$

Inequality (2.22) guarantees that the functional $V: L^\infty([-r,0);\Re^n) \to \Re_+$ is Lipschitz on bounded sets of the state space $L^\infty([-r,0);\Re^n)$. However, inequality (2.22) does not guarantee Frechet differentiability of the functional $V: L^\infty([-r,0);\Re^n) \to \Re_+$ nor that the limit $\lim_{t \to 0^+} t^{-1}(V(x_t) - V(x_0))$ exists for the solution $x(t)$ of (1.1) with initial condition $x(t) = x_0(t)$ for $t \in [-r,0)$ a.e., corresponding to input $d \in L^\infty([-r,+\infty);D)$. Notice that inequality (2.23) guarantees that for every $x_0 \in L^\infty([-r,0);\Re^n)$, $d \in L^\infty([-r,+\infty);D)$, the solution $x(t)$ of (1.1) with initial condition $x(t) = x_0(t)$ for $t \in [-r,0)$ a.e., corresponding to input $d \in L^\infty([-r,+\infty);D)$ satisfies $\limsup_{t \to 0^+} t^{-1}(V(x_t) - V(x_0)) \leq -V(x_0)$.

Property 4: Sufficient Conditions for Stability Properties

Theorem 2.5 is not the most convenient way of proving RGAS or ISS for (1.1). For practical purposes we can use the following result, which is an extension of the classical Razumikhin theorem for time delay systems (see [8,10]).

**Theorem 2.6:** *Consider system (1.1) under assumptions (H1), (H2), (H3). Assume that there exists a continuous, positive definite and radially unbounded function $W: \Re^n \to \Re_+$, a continuous, non-decreasing function $\gamma: \Re_+ \to \Re_+$ and a constant $\lambda \in (0,1)$ such that the following inequality holds for all $x \in L^\infty([-r,0);\Re^n)$, $d \in L^\infty([-r,0];D)$, $u \in L^\infty([-r,0];U)$:*

$$W(f(x,d,u)) \leq \lambda \sup_{-r \leq s < 0} W(x(s)) + \gamma(\|u\|) \qquad (2.24)$$

*Then system (1.1) is ISS from the input $u \in U \subseteq \Re^{m_2}$ uniformly in $d \in D \subseteq \Re^{m_1}$.*

**Proof:** Define the functional:

$$V(x) := \sup_{-r \leq s < 0} \exp(\sigma s) W(x(s)), \text{ for all } x \in L^\infty([-r,0);\Re^n) \qquad (2.25)$$

where $\sigma \in \left(0, \frac{1}{r}\ln\left(\frac{1}{\lambda}\right)\right)$ is a constant. Since $W: \Re^n \to \Re_+$ is a continuous, positive definite and radially unbounded function, there exist functions $a_1, a_2 \in K_\infty$ such that:

$$a_1(|x|) \leq W(x) \leq a_2(|x|), \text{ for all } x \in \Re^n \qquad (2.26)$$

Using (2.25) and (2.26) we conclude that the following inequality holds:

$$\tilde{a}_1(\|x\|) \leq V(x) \leq a_2(\|x\|), \text{ for all } x \in L^\infty([-r,0);\Re^n) \qquad (2.27)$$

where $\tilde{a}_1(s) := \exp(-\sigma r) a_1(s)$, for all $s \geq 0$.

Let arbitrary $x_0 \in L^\infty([-r,0);\Re^n)$, $d \in L^\infty([-r,+\infty);D)$, $u \in L^\infty([-r,+\infty);U)$ and consider the unique solution $x(t)$ of (1.1) with initial condition $x(t) = x_0(t)$ for $t \in [-r,0)$ a.e., corresponding to inputs $d \in L^\infty([-r,+\infty);D)$, $u \in L^\infty_{loc}([-r,+\infty);U)$. Let arbitrary $t \in (0, t_{\max}(x_0,d,u))$, where $t_{\max}(x_0,d,u) > 0$ is the maximal existence time of the solution and let $h \in (0, t_{\max}(x_0,d,u) - t)$. Using definition (2.25) we get:

$$
\begin{aligned}
V(x_{t+h}) &:= \sup_{-r \leq s < 0} \exp(\sigma s) W(x(t+h+s)) \\
&= \max\left( \sup_{-r \leq s \leq -h} \exp(\sigma s) W(x(t+h+s)), \sup_{-h \leq s < 0} \exp(\sigma s) W(x(t+h+s)) \right) \\
&= \max\left( \exp(-\sigma h) \sup_{-r+h \leq s < 0} \exp(\sigma s) W(x(t+s)), \sup_{-h \leq s < 0} \exp(\sigma s) W(x(t+h+s)) \right) \\
&\leq \max\left( \exp(-\sigma h) V(x_t), \sup_{0 \leq s < h} W(x(t+s)) \right)
\end{aligned}
\quad (2.28)
$$

Using (2.24) we get for almost all $q \in [0,h)$:

$$
W(x(t+q)) = W\big(f(x_{t+q}, d_{t+q}, u_{t+q})\big) \leq \lambda \max\left( \sup_{-r+q \leq s < 0} W(x(t+s)), \sup_{0 \leq s < q} W(x(t+s)) \right) + \gamma(\|u\|) \quad (2.29)
$$

where $\|u\| = \sup_{t \geq 0} \|u_t\|$. Definition (2.25) and inequality (2.29) gives:

$$
\sup_{0 \leq q < h} W(x(t+q)) \leq \lambda \max\left( \exp(\sigma r) V(x_t), \sup_{0 \leq s < h} W(x(t+s)) \right) + \gamma(\|u\|)
$$

which directly implies:

$$
\sup_{0 \leq q < h} W(x(t+q)) \leq \lambda \exp(\sigma r) V(x_t) + \frac{1}{1-\lambda} \gamma(\|u\|) \quad (2.30)
$$

Since $\lambda \exp(\sigma r) \leq \exp(-\sigma h)$ for all $h \leq \frac{1}{\sigma} \ln\left(\frac{1}{\lambda}\right) - r$ (notice that $\sigma \in \left(0, \frac{1}{r} \ln\left(\frac{1}{\lambda}\right)\right]$), we get from (2.28), (2.30) for all $h \in (0, t_{\max}(x_0, d, u) - t)$ with $h \leq \frac{1}{\sigma} \ln\left(\frac{1}{\lambda}\right) - r$:

$$
V(x_{t+h}) \leq \exp(-\sigma h) V(x_t) + \frac{1}{1-\lambda} \gamma(\|u\|) \quad (2.31)
$$

Using the Boundedness-Implies-Continuation property, inequalities (2.27), (2.31) and a standard contradiction argument, we conclude that $t_{\max}(x_0, d, u) = +\infty$. Therefore, we conclude that (2.31) holds for all $t \geq 0$ and $h \geq 0$ with $h \leq \frac{1}{\sigma} \ln\left(\frac{1}{\lambda}\right) - r$.

Using induction and (2.31) we show that $V(x_{ih}) \leq \exp(-i\sigma h) V(x_0) + \frac{1}{1-\lambda} \gamma(\|u\|) \sum_{j=0}^{i-1} \exp(-j\sigma h)$, for all integers $i \geq 1$ and $h = \frac{1}{\sigma} \ln\left(\frac{1}{\lambda}\right) - r$. Consequently, the previous inequality in conjunction with (2.31) shows that $V(x_{ih+q}) \leq \exp(-\sigma(q+ih)) V(x_0) + \frac{1}{1-\lambda} \gamma(\|u\|) \frac{2 - \exp(-\sigma h)}{1 - \exp(-\sigma h)}$, for all integers $i \geq 0$, $h = \frac{1}{\sigma} \ln\left(\frac{1}{\lambda}\right) - r$ and for all $q \in \left[0, \frac{1}{\sigma} \ln\left(\frac{1}{\lambda}\right) - r\right)$. Since, for every $t \geq 0$ there exists integer $i \geq 0$ and $q \in \left[0, \frac{1}{\sigma} \ln\left(\frac{1}{\lambda}\right) - r\right)$ with $t = ih + q$, we obtain:

$$
V(x_t) \leq \exp(-\sigma t) V(x_0) + \frac{1}{1-\lambda} \frac{2 - \lambda \exp(\sigma r)}{1 - \lambda \exp(\sigma r)} \gamma(\|u\|), \text{ for all } t \geq 0 \quad (2.32)
$$

Inequality (2.32) in conjunction with (2.27) gives for all $\varepsilon > 0$:

$$\tilde{a}_1(\|x_t\|) \le \max\left((1+\varepsilon)\exp(-\sigma t)a_2(\|x_0\|), \frac{1+\varepsilon^{-1}}{1-\lambda}\frac{2-\lambda\exp(\sigma r)}{1-\lambda\exp(\sigma r)}\gamma(\|u\|)\right), \text{ for all } t \ge 0 \quad (2.32)$$

Inequality (2.32) shows that there exists $\sigma \in KL$ and a continuous, non-decreasing function $\tilde{\gamma}: \Re_+ \to \Re_+$ such that for every $x_0 \in L^\infty([-r,0]; \Re^n)$, $d \in L^\infty([-r,+\infty); D)$, $u \in L^\infty([-r,+\infty); U)$, there exists $x \in L^\infty_{loc}([-r,+\infty); \Re^n)$ satisfying $x(t) = f(x_t, d_t, u_t)$ for $t \in [0,+\infty)$ a.e., $x(t) = x_0(t)$ for $t \in [-r,0)$ a.e. and $\|x_t\| \le \sigma(\|x_0\|,t) + \tilde{\gamma}(\|u\|)$ for all $t \ge 0$. A standard causality argument shows that the estimate $\|x_t\| \le \sigma(\|x_0\|,t) + \sup_{0\le s \le t}\tilde{\gamma}(\|u_s\|)$ holds for all $t \ge 0$ and $u \in L^\infty_{loc}([-r,+\infty); U)$. Therefore, Theorem 2.4 implies that system (1.1) is ISS from the input $u \in U \subseteq \Re^{m_2}$ uniformly in $d \in D \subseteq \Re^{m_1}$. The proof is complete. ◁

**Example 2.7:** Consider the linear system described by the single IDE:

$$x(t) = d(t)\int_{-1}^{0} q(s)x(t+s)ds + u(t) \quad (2.33)$$
$$x(t) \in \Re, d(t) \in [-1,1], u(t) \in \Re$$

where $q:[-1,0] \to \Re$ is a continuous function. Using Theorem 2.6 with $W(x) = |x|$, we can conclude that system (2.33) is ISS from the input $u \in \Re$ under the assumption that there exists $\lambda \in (0,1)$ such that $\int_{-1}^{0} |q(s)|ds \le \lambda$. In fact, following the proof of Theorem 2.6, we can show that the gain function of the input $u \in \Re$ is linear. ◁

# 3. Systems Described by a Single First-Order Hyperbolic Partial Differential Equation

Consider the system described by a single FOH-PDE of the form (1.2), (1.3), (1.4) under the following assumption:

**(A1)** *There exist non-decreasing functions $L: \Re_+ \to \Re_+$, $\sigma: \Re_+ \to \Re_+$ such that the functionals $G: \Re^m \times L^\infty((0,1]; \Re) \to \Re$, $K_i: \Re^m \times L^\infty((0,1]; \Re) \to \Re$ ($i=1,...,N$) satisfy the following inequalities for all $R \ge 0$:*

$$|G(w,x) - G(w,y)| + \sum_{i=1}^{N}|K_i(w,x) - K_i(w,y)| \le L(R)h \sup_{0<z\le h}(|x(z)-y(z)|) + L(R)\sup_{h<z\le 1}(|x(z)-y(z)|),$$
$$\text{for all } h \in (0,1), w \in W \text{ with } |w| \le R \text{ and } x,y \in L^\infty((0,1];[-R,R]) \quad (3.1)$$

$$|G(w,x)| + \sum_{i=1}^{N}|K_i(w,x)| \le \sigma(R),$$
$$\text{for all } w \in W \text{ with } |w| \le R \text{ and } x \in L^\infty((0,1];[-R,R]) \quad (3.2)$$

We show next how to transform system (1.2), (1.3), (1.4) under assumption (A1) to an equivalent system described by IDEs.

Consider a classical solution of (1.2), (1.3), (1.4) with sufficiently regular initial condition $x(0,z) = x_0(z)$ for $z \in (0,1]$ corresponding to a sufficiently regular input $w: \Re_+ \to W$. We assume that the solution exists for all $t \in J$, where $J \subseteq \Re_+$ is an interval with non-empty interior and $0 \in J$. First, we define for $t \in J$:

$$v(t) = G(w(t), x_t) \tag{3.3}$$

If we define:

$$v(-q) = \exp\left(-\int_{-q}^{0} a(p(s), c(q+s))ds\right) x_0(cq) - \int_{-q}^{0} \exp\left(-\int_{-q}^{\tau} a(p(s), c(q+s))ds\right) g(c(q+\tau)) p(\tau) d\tau,$$
$$\text{for almost all } q \in (0, c^{-1}] \tag{3.4}$$

then integration on the characteristic line allows us to conclude that the solution of (1.2), (1.3), (1.4) must satisfy the following equation for all $t \in J$ and $z \in [0,1]$:

$$x(t,z) = \exp\left(\int_{t-c^{-1}z}^{t} a(p(s), z+c(s-t))ds\right) v(t-c^{-1}z) + \int_{t-c^{-1}z}^{t} \exp\left(\int_{\tau}^{t} a(p(s), z+c(s-t))ds\right) g(z+c(\tau-t)) p(\tau) d\tau \tag{3.5}$$

Let $p_t \in L^\infty\left([-c^{-1}, 0); \Re^N\right)$ be the mapping defined by $(p_t)(s) = p(t+s)$ for $t \in J$, $s \in [-c^{-1}, 0)$. Similarly, let $v_t \in L^\infty\left([-c^{-1}, 0); \Re\right)$ be the mapping defined by $(v_t)(s) = v(t+s)$ for $t \in J$, $s \in [-c^{-1}, 0)$. Clearly, equation (3.5) allows us to define for every $p \in L^\infty\left([-c^{-1}, 0); \Re^N\right)$ a bounded linear operator $A(p): L^\infty\left([-c^{-1}, 0); \Re\right) \to L^\infty((0,1]; \Re)$ and a mapping $B: L^\infty\left([-c^{-1}, 0); \Re^N\right) \to L^\infty((0,1]; \Re)$, so that the following equation holds for all $t \in J$:

$$x_t = A(p_t) v_t + B(p_t) \tag{3.6a}$$

where

$$(A(p)v)(z) := \exp\left(\int_{-c^{-1}z}^{0} a(p(t+q), z+cq)dq\right) v(-c^{-1}z)$$

$$(B(p))(z) := \int_{-c^{-1}z}^{0} \exp\left(\int_{s}^{0} a(p(q), z+cq)dq\right) g(z+cs) p(s) ds$$
, for $z \in (0,1]$ \hfill (3.6b)

It follows from (3.3) and (3.6a) that the following equations hold for all $t \in J$:

$$p_i(t) = F_i(w(t), p_t, v_t), \text{ for } i = 1, \ldots, N \text{ and } v(t) = F_{N+1}(w(t), p_t, v_t) \tag{3.7}$$

where the mappings $F_i: W \times L^\infty\left([-c^{-1}, 0); \Re^{N+1}\right) \to \Re$ $(i = 1, \ldots, N+1)$ are defined for every $(w, p, v) \in W \times L^\infty\left([-c^{-1}, 0); \Re^{N+1}\right)$ by the equations:

$$F_i(w, p, v) = K_i(w, A(p)v + B(p)), \text{ for } i = 1, \ldots, N \text{ and } F_{N+1}(w, p, v) = G(w, A(p)v + B(p)) \tag{3.8}$$

In other words, every classical solution $x_t$ of (1.2), (1.3), (1.4) corresponds to a solution $(p_t, v_t) \in L^\infty\left([-c^{-1}, 0); \Re^{N+1}\right)$ of the system described by the integral delay equations (3.7), (3.8).

The procedure described above is rigorous and the system described by the IDEs (3.7), (3.8) is equivalent to the system described by the single FOH-PDE (1.2), (1.3), (1.4). This is shown by the following theorem.

**Theorem 3.1:** *Assume that $a \in C^0(\Re^N \times [0,1]; \Re)$ is locally Lipschitz and $g_i \in C^0([0,1]; \Re)$ ($i = 1,...,N$). Moreover, assume that there exist non-decreasing functions $L: \Re_+ \to \Re_+$, $\sigma: \Re_+ \to \Re_+$ such that the functionals $G: W \times L^\infty((0,1]; \Re) \to \Re$, $K_i: W \times L^\infty((0,1]; \Re) \to \Re$ ($i = 1,...,N$) satisfy inequalities (3.1), (3.2) for all $R \geq 0$. Finally, assume that the following assumption holds for the mapping $F: W \times L^\infty\left([-c^{-1}, 0); \Re^{N+1}\right) \to \Re^{N+1}$, where $F(w, p, v) := (F_1(w, p, v), ..., F_{N+1}(w, p, v))$ and $F_i: W \times L^\infty\left([-c^{-1}, 0); \Re^{N+1}\right) \to \Re$ ($i = 1,...,N+1$) are defined by (3.8) with the aid of (3.6b).*

**(A2)** *For every $\delta > 0$, $(p, v) \in L^\infty\left([-c^{-1}, \delta); \Re^{N+1}\right)$, $w \in L^\infty([0, \delta]; W)$, every function $\xi: [-c^{-1}, \delta) \to \Re^{N+1}$ defined by $\xi(t) = (p(t), v(t))$ for $t \in [-c^{-1}, 0)$ a.e. and $\xi(t) = F(w(t), p_t, v_t)$ for $t \in [0, \delta)$ a.e. satisfies $\xi \in L^\infty\left([-c^{-1}, \delta); \Re^{N+1}\right)$.*

*Then for every $x_0 \in L^\infty([0,1]; \Re)$, $w \in L^\infty_{loc}([0, +\infty); W)$ there exists $t_{\max} = t_{\max}(x_0, w) \in (0, +\infty]$ and a unique locally bounded mapping $[0, t_{\max}) \ni t \to x_t \in L^\infty([0,1]; \Re)$ satisfying (1.3), (1.4), (3.3) for $t \in [0, t_{\max})$ a.e. with $(p, v) \in L^\infty_{loc}([0, t_{\max}); \Re^{N+1})$, and*

$$x(t, z) = \int_{\max(0, t - c^{-1}z)}^{t} \exp\left(\int_\tau^t a(p(s), z + c(s-t)) ds\right) g(z + c(\tau - t)) p(\tau) d\tau$$
$$+ \exp\left(\int_{\max(0, t - c^{-1}z)}^{t} a(p(s), z + c(s-t)) ds\right) \begin{cases} v(t - c^{-1}z) & \text{for } z \leq ct \\ x_0(z - ct) & \text{for } z > ct \end{cases} \qquad (3.9)$$

*for all $t \in [0, t_{\max})$ and almost all $z \in [0,1]$. The mapping $t \to x_t \in L^\infty([0,1]; \Re)$ is obtained from (3.5) and the solution of the system (3.7) with arbitrary initial condition $(p, v) \in L^\infty\left([-c^{-1}, 0); \Re^{N+1}\right)$ that satisfies (3.4). Moreover, the mapping $t \to x_t \in L^\mu([0,1]; \Re)$ is continuous on $[0, t_{\max})$ for every $\mu \in [1, +\infty)$ and if $t_{\max} < +\infty$ then $\limsup_{t \to t_{\max}^-} \sup_{0 \leq z \leq 1} |(x_t)(z)| = +\infty$. Finally, there exist non-decreasing functions $P: \Re_+ \to \Re_+$, $Q: \Re_+ \to \Re_+$ such that for every $x_0 \in L^\infty([0,1]; \Re)$, $y_0 \in L^\infty([0,1]; \Re)$, $w \in L^\infty_{loc}([0, +\infty); W)$ it holds that:*

$$\sup_{0 \leq z \leq 1} |(x_t)(z) - (y_t)(z)| \leq Q(s(t)) \exp(P(s(t))t) \sup_{0 \leq z \leq 1} |(x_0)(z) - (y_0)(z)|, \text{ for all } t \in [0, \delta) \qquad (3.10)$$

*where $\delta = \min(t_{\max}(x_0, w), t_{\max}(y_0, w))$, $s(t) := \max\left(\sup_{0 \leq \tau < t} \sup_{0 \leq z \leq 1} |(x_\tau)(z)|, \sup_{0 \leq \tau < t} \sup_{0 \leq z \leq 1} |(y_\tau)(z)|, \sup_{0 \leq \tau \leq t} |w(t)|\right)$ and the mappings $t \to x_t \in L^\infty([0,1]; \Re)$, $t \to y_t \in L^\infty([0,1]; \Re)$ are the unique mappings satisfying (1.3), (1.4), (3.3) for $t \in [0, t_{\max})$ a.e. with $(p, v) \in L^\infty_{loc}([0, t_{\max}); \Re^{N+1})$, (3.9) for all $t \in [0, t_{\max})$ and almost all $z \in [0,1]$ as well as the analogous equations with $y$ replacing $x$.*

The reader should notice that no regularity result is provided for the mapping $t \to x_t \in L^\infty([0,1];\Re)$. It is clear that the unique mapping $t \to x_t \in L^\infty([0,1];\Re)$ satisfying (1.3), (1.4), (3.3) for $t \in [0,t_{\max})$ a.e. with $(p,v) \in L^\infty_{loc}([0,t_{\max});\Re^{N+1})$, (3.9) for all $t \in [0,t_{\max})$ and almost all $z \in [0,1]$ is not a classical solution of the initial value problem (1.2), (1.3), (1.4) with initial condition $x(0,z) = x_0(z)$ for $z \in [0,1]$ a.e., where $x_0 \in L^\infty([0,1];\Re)$. The mapping $t \to x_t \in L^\infty([0,1];\Re)$ satisfying (1.3), (1.4), (3.3) for $t \in [0,t_{\max})$ a.e. with $(p,v) \in L^\infty_{loc}([0,t_{\max});\Re^{N+1})$, (3.9) for all $t \in [0,t_{\max})$ and almost all $z \in [0,1]$ is close to the notion of the "mild" solution for first-order systems of hyperbolic PDEs (see [7] and references therein).

**Proof:** The proof of the theorem is long, incorporating several distinct results, so we divide it in four steps.

Step 1: We show that system (3.7) satisfies assumptions (H1), (H2). Therefore, by virtue of Theorem 2.1 for every input $w \in L^\infty_{loc}([0,+\infty);W)$ and for every initial condition $(p,v) \in L^\infty([-c^{-1},0);\Re^{N+1})$ that satisfies (3.4) for certain $x_0 \in L^\infty([0,1];\Re)$, system (3.7) admits a unique solution. The obtained solution of system (3.7) gives us directly (by means of (3.4) and (3.5)) a locally bounded mapping $[0,t_{\max}) \ni t \to x_t \in L^\infty([0,1];\Re)$ satisfying (1.3), (1.4), (3.3) for $t \in [0,t_{\max})$ a.e. with $(p,v) \in L^\infty_{loc}([0,t_{\max});\Re^{N+1})$, and (3.9) for all $t \in [0,t_{\max})$ and almost all $z \in [0,1]$.

Step 2: We show that for every $w \in L^\infty_{loc}([0,+\infty);W)$, $x_0 \in L^\infty([0,1];\Re)$ there exists at most one locally bounded mapping $[0,t_{\max}) \ni t \to x_t \in L^\infty([0,1];\Re)$ satisfying (1.3), (1.4), (3.3) for $t \in [0,t_{\max})$ a.e. with $(p,v) \in L^\infty_{loc}([0,t_{\max});\Re^{N+1})$, and (3.9) for all $t \in [0,t_{\max})$ and almost all $z \in [0,1]$. Therefore, the mapping obtained in Step 1 is the unique mapping that satisfies (1.3), (1.4), (3.3) for $t \in [0,t_{\max})$ a.e. with $(p,v) \in L^\infty_{loc}([0,t_{\max});\Re^{N+1})$, and (3.9) for all $t \in [0,t_{\max})$ and almost all $z \in [0,1]$. This means that for every initial condition $(p,v) \in L^\infty([-c^{-1},0);\Re^{N+1})$ that satisfies (3.4) for the given $x_0 \in L^\infty([0,1];\Re)$ the solution of system (3.7) will be the same.

Step 3: Making use of the results of Theorem 2.1, we show that if $t_{\max} < +\infty$ then $\limsup_{t \to t_{\max}^-} \sup_{0 \le z \le 1} |(x_t)(z)| = +\infty$ and that there exist non-decreasing functions $P: \Re_+ \to \Re_+$, $Q: \Re_+ \to \Re_+$ such that for every $x_0 \in L^\infty([0,1];\Re)$, $y_0 \in L^\infty([0,1];\Re)$, $w \in L^\infty_{loc}([0,+\infty);W)$ estimate (3.10) holds.

Step 4: We show that for every $w \in L^\infty_{loc}([0,+\infty);W)$, $x_0 \in L^\infty([0,1];\Re)$ the locally bounded mapping $[0,t_{\max}) \ni t \to x_t \in L^\infty([0,1];\Re)$ satisfying (1.3), (1.4), (3.3) for $t \in [0,t_{\max})$ a.e. with $(p,v) \in L^\infty_{loc}([0,t_{\max});\Re^{N+1})$, and (3.9) for all $t \in [0,t_{\max})$ and almost all $z \in [0,1]$ gives us a mapping $t \to x_t \in L^\mu([0,1];\Re)$ which is continuous on $[0,t_{\max})$ for every $\mu \in [1,+\infty)$.

Step 1: Assumption (A2) guarantees that assumption (H2) holds for system (3.7). Likewise, the validity of assumption (H1) for system (3.7) is a direct consequence of assumption (A1). To see this, let $R \ge 0$, $(p,v) \in L^\infty([-c^{-1},0);\Re^{N+1})$, $(\tilde{p},\tilde{v}) \in L^\infty([-c^{-1},0);\Re^{N+1})$ with $\|(p,v)\| \le R$, $\|(\tilde{p},\tilde{v})\| \le R$ be given and let $x,\tilde{x} \in L^\infty((0,1];\Re)$ be defined by $x = A(p)v + B(p)$, $\tilde{x} = A(\tilde{p})\tilde{v} + B(\tilde{p})$, where $A(p): L^\infty([-c^{-1},0);\Re) \to L^\infty((0,1];\Re)$ and $B: L^\infty([-c^{-1},0);\Re^N) \to L^\infty((0,1];\Re)$ are the mappings defined by (3.6b), i.e.,

$$x(z) = \exp\left(\int_{-c^{-1}z}^{0} a(p(s), z+cs)ds\right) v(-c^{-1}z) + \int_{-c^{-1}z}^{0} \exp\left(\int_{\tau}^{0} a(p(s), z+cs)ds\right) g(z+c\tau)p(\tau)d\tau, \text{ for } z \in [0,1] \text{ a.e. } (3.11)$$

$$\widetilde{x}(z) = \exp\left(\int_{-c^{-1}z}^{0} a(\widetilde{p}(s), z+cs)ds\right) \widetilde{v}(-c^{-1}z) + \int_{-c^{-1}z}^{0} \exp\left(\int_{\tau}^{0} a(\widetilde{p}(s), z+cs)ds\right) g(z+c\tau)\widetilde{p}(\tau)d\tau, \text{ for } z \in [0,1] \text{ a.e. } (3.12)$$

Define:

$$\widetilde{M}(R) := \max_{z \in [0,1], |p| \leq R} |a(p,z)| \text{ and } C := \max_{z \in [0,1]} |g(z)| \quad (3.13)$$

and notice that (3.11) and definitions (3.13) imply:

$$\max\left(\sup_{0 \leq z \leq 1} |x(z)|, \sup_{0 \leq z \leq 1} |\widetilde{x}(z)|\right) \leq (1 + c^{-1}C) \exp(c^{-1}\widetilde{M}(R))R \quad (3.14)$$

Using (3.2), definitions (3.8) and (3.14), we obtain:

$$|F(w, p, v)| \leq \sigma\left((1 + c^{-1}C) \exp(c^{-1}\widetilde{M}(R))R\right),$$
$$\text{for all } (w, p, v) \in W \times L^{\infty}([-c^{-1}, 0); \Re^{N+1}) \text{ with } \|(p,v)\| \leq R \text{ and } |w| \leq R \quad (3.15)$$

Inequality (3.15) shows that inequality (2.2) holds for the non-decreasing function $a(R) := \sigma\left((1 + c^{-1}C) \exp(c^{-1}\widetilde{M}(R))R\right)$ for all $R \geq 0$ (notice that definition (3.13) implies that $\widetilde{M}(R)$ is non-decreasing).

Using (3.1), definitions (3.8) and (3.14), we obtain:

$$|F(w, p, v) - F(w, \widetilde{p}, \widetilde{v})| \leq L(K(R))h \sup_{0 < z \leq h}(|x(z) - \widetilde{x}(z)|) + L(K(R)) \sup_{h < z \leq 1}(|x(z) - \widetilde{x}(z)|),$$
$$\text{for all } h \in (0,1), w \in W \text{ with } |w| \leq R \quad (3.16)$$

where $K(R) := (1 + c^{-1}C) \exp(c^{-1}\widetilde{M}(R))R$. Define $L_a(R) := \sup\left\{\frac{|a(p,z) - a(\widetilde{p},z)|}{|p - \widetilde{p}|} : p \neq \widetilde{p}, |p| \leq R, |\widetilde{p}| \leq R, z \in [0,1]\right\}$ and notice that $L_a(R) < +\infty$ for all $R \geq 0$ (since $a \in C^0(\Re^N \times [0,1]; \Re)$ is assumed to be locally Lipschitz). Using (3.11), (3.12) and definitions (3.13) we obtain for $z \in [0,1]$ a.e.:

$$|x(z) - \widetilde{x}(z)| \leq B_1(R)|v(-c^{-1}z) - \widetilde{v}(-c^{-1}z)| + B_2(R) \int_{-c^{-1}z}^{0} |p(s) - \widetilde{p}(s)|ds \quad (3.17)$$

where $B_1(R) := \exp(\widetilde{M}(R)c^{-1})$, $B_2(R) := \left(RL_a(R) + C + CRc^{-1}L_a(R)\right)\exp(3\widetilde{M}(R)c^{-1})$. Combining (3.16) and (3.17) we obtain for all $h \in (0,1)$, $w \in W$ with $|w| \leq R$:

$$\begin{aligned}
|F(w,p,v) - F(w,\tilde{p},\tilde{v})| \leq\;& L(K(R))B_1(R)h \sup_{-c^{-1}h \leq s < 0} (|v(s) - \tilde{v}(s)|) \\
& + 2L(K(R))B_2(R)c^{-1}h \sup_{-c^{-1}h \leq s < 0} |p(s) - \tilde{p}(s)| \\
& + L(K(R))B_1(R) \sup_{c^{-1} \leq s < -c^{-1}h} (|v(s) - \tilde{v}(s)|) \\
& + L(K(R))B_2(R)c^{-1} \sup_{c^{-1} \leq s < -c^{-1}h} |p(s) - \tilde{p}(s)|
\end{aligned} \quad (3.18)$$

Inequality (3.18) shows that inequality (2.1) holds for appropriate non-decreasing functions $M:\Re_+ \to \Re_+$, $N:\Re_+ \to \Re_+$.

Let arbitrary $x_0 \in L^\infty([0,1];\Re)$ and $w \in L^\infty_{loc}([0,+\infty);W)$. Consider now the solution of system (3.7) corresponding to input $w \in L^\infty_{loc}([0,+\infty);W)$ with initial condition $(p,v) \in L^\infty([-c^{-1},0);\Re^{N+1})$ with $p(s) \equiv 0$ and (3.4), i.e.,

$$v(t) = \exp\left(-\int_t^0 a(0, c(s-t))ds\right) x_0(-ct), \text{ for almost all } t \in [-c^{-1}, 0) \quad (3.19)$$

Theorem 2.1 guarantees that there exists $t_{\max} := t_{\max}(x_0, w) \in (0, +\infty]$ so that system (3.7) admits a unique solution $(p,v) \in L^\infty_{loc}([-c^{-1}, t_{\max});\Re^{N+1})$. Using (3.5) and (3.19), we obtain a locally bounded mapping $[0, t_{\max}) \ni t \to x_t \in L^\infty([0,1];\Re)$, which satisfies (1.3), (1.4), (3.3) for $t \in [0, t_{\max})$ a.e. and (3.9) for all $t \in [0, t_{\max})$ and almost all $z \in [0,1]$.

Step 2: We show that by contradiction that for every $w \in L^\infty_{loc}([0,+\infty);W)$, $x_0 \in L^\infty([0,1];\Re)$ there exists at most one locally bounded mapping $[0, t_{\max}) \ni t \to x_t \in L^\infty([0,1];\Re)$ satisfying (1.3), (1.4), (3.3) for $t \in [0, t_{\max})$ a.e. with $(p,v) \in L^\infty_{loc}([0, t_{\max});\Re^{N+1})$, and (3.9) for all $t \in [0, t_{\max})$ and almost all $z \in [0,1]$.

Suppose that there exist $w \in L^\infty_{loc}([0,+\infty);W)$, $x_0 \in L^\infty([0,1];\Re)$, a constant $\delta > 0$ and two locally bounded mappings $[0, \delta) \ni t \to x_t \in L^\infty([0,1];\Re)$, $[0, \delta) \ni t \to \tilde{x}_t \in L^\infty([0,1];\Re)$ so that (3.9) holds for all $t \in [0, \delta)$ and almost all $z \in [0,1]$, equation

$$\begin{aligned}
\tilde{x}(t,z) =\;& \int_{\max(0, t-c^{-1}z)}^{t} \exp\left(\int_\tau^t a(\tilde{p}(s), z + c(s-t))ds\right) g(z + c(\tau-t))\tilde{p}(\tau)d\tau \\
& + \exp\left(\int_{\max(0, t-c^{-1}z)}^{t} a(\tilde{p}(s), z + c(s-t))ds\right) \begin{cases} \tilde{v}(t - c^{-1}z) & \text{for } z \leq ct \\ x_0(z - ct) & \text{for } z > ct \end{cases}
\end{aligned} \quad (3.20)$$

for all $t \in [0, \delta)$ and almost all $z \in [0,1]$, (1.3), (1.4), (3.3) hold for $t \in [0, \delta)$ a.e., equations

$$\tilde{p}_i(t) = K_i(w(t), \tilde{x}_t), \text{ for } i=1,\ldots,N \text{ and } \tilde{p}(t) = (\tilde{p}_1(t), \ldots, \tilde{p}_N(t))' \in \Re^N \quad (3.21)$$

$$\tilde{v}(t) = G(w(t), \tilde{x}_t) \quad (3.22)$$

hold for $t \in [0, \delta)$ a.e. and such that

$$\sup_{t\in[0,\delta)}\sup_{z\in[0,1]}\left|(x_t)(z)-(\widetilde{x}_t)(z)\right|>0 \qquad (3.23)$$

Let $t_1\in[0,\delta)$ be the greatest lower bound of all times $t\in[0,\delta)$ for which $\sup_{z\in[0,1]}\left|(x_t)(z)-(\widetilde{x}_t)(z)\right|>0$. Notice that by virtue of (3.23), $t_1\in[0,\delta)$ is well defined and that $\sup_{t\in[0,t_1)}\sup_{z\in[0,1]}\left|(x_t)(z)-(\widetilde{x}_t)(z)\right|=0$ for the case that $t_1>0$.

We show next that there exists $h>0$ so that $\sup_{t\in[0,t_1+h)}\sup_{z\in[0,1]}\left|(x_t)(z)-(\widetilde{x}_t)(z)\right|=0$, which contradicts the fact that $t_1\in[0,\delta)$ is the greatest lower bound of all times $t\in[0,\delta)$ for which $\sup_{z\in[0,1]}\left|(x_t)(z)-(\widetilde{x}_t)(z)\right|>0$.

An important property that is used at this point is the so-called classical semigroup property (see [10]), i.e., the fact that every locally bounded mapping $[0,t_{\max})\ni t\to x_t\in L^\infty([0,1];\Re)$ satisfying (1.3), (1.4), (3.3) for $t\in[0,t_{\max})$ a.e. with $(p,v)\in L^\infty_{loc}\left([0,t_{\max});\Re^{N+1}\right)$, and (3.9) for all $t\in[0,t_{\max})$ and almost all $z\in[0,1]$ satisfies the following identity for all $t\in[0,t_{\max})$ and for all $t_0\in[0,t]$:

$$x(t,z) = \int_{\max(0,t-t_0-c^{-1}z)}^{t-t_0} \exp\left(\int_\tau^{t-t_0} a(p(t_0+s),z+c(s+t_0-t))ds\right) g(z+c(\tau+t_0-t))p(t_0+\tau)d\tau$$
$$+\exp\left(\int_{\max(0,t-t_0-c^{-1}z)}^{t-t_0} a(p(t_0+s),z+c(s+t_0-t))ds\right)\begin{cases} v(t-c^{-1}z) & \text{for } z+ct_0\leq ct \\ (x_{t_0})(z+ct_0-ct) & \text{for } z+ct_0>ct \end{cases},$$
$$\text{for } z\in[0,1] \text{ a.e.} \qquad (3.24)$$

Therefore, identity (3.24) holds for all $t\in[0,\delta)$ and for all $t_0\in[0,t]$. Similarly, the following identity

$$\widetilde{x}(t,z) = \int_{\max(0,t-t_0-c^{-1}z)}^{t-t_0} \exp\left(\int_\tau^{t-t_0} a(\widetilde{p}(t_0+s),z+c(s+t_0-t))ds\right) g(z+c(\tau+t_0-t))\widetilde{p}(t_0+\tau)d\tau$$
$$+\exp\left(\int_{\max(0,t-t_0-c^{-1}z)}^{t-t_0} a(\widetilde{p}(t_0+s),z+c(s+t_0-t))ds\right)\begin{cases} \widetilde{v}(t-c^{-1}z) & \text{for } z+ct_0\leq ct \\ (\widetilde{x}_{t_0})(z+ct_0-ct) & \text{for } z+ct_0>ct \end{cases},$$
$$\text{for } z\in[0,1] \text{ a.e.} \qquad (3.25)$$

holds for all $t\in[0,\delta)$ and for all $t_0\in[0,t]$. Using (3.24) and (3.25) with $t_0\in(0,t_1)$ for the case $t_1>0$ and $t_0=0$ for the case $t_1=0$ in conjunction with (1.3), (3.3), (3.13), (3.21), (3.22) (which implies that $p(s)=\widetilde{p}(s)$, $v(s)=\widetilde{v}(s)$, $x_s=\widetilde{x}_s$ for all $s\in[0,t_0]$), we obtain:

$$\left|x(t,z)-\widetilde{x}(t,z)\right|\leq$$
$$(t-t_0)\exp\left(3\widetilde{M}(R)(t-t_0)\right)\left(C+CRL_a(R)(t-t_0)+2RL_a(R)\right)\sup_{t_0\leq s\leq t}\left|p(s)-\widetilde{p}(s)\right| \qquad (3.26)$$
$$+\exp\left(\widetilde{M}(R)(t-t_0)\right)\begin{cases} \left|v(t-c^{-1}z)-\widetilde{v}(t-c^{-1}z)\right| & \text{for } z+ct_0\leq ct \\ 0 & \text{for } z+ct_0>ct \end{cases}$$

for $z \in [0,1]$ a.e. and for all $t \in \left[t_0, \frac{t_1+\delta}{2}\right]$, where

$$R = \sup_{0 \leq s \leq \frac{t_1+\delta}{2}} \max\left(|p(s)|, |\tilde{p}(s)|, |v(s)|, |\tilde{v}(s)|, \sup_{0 \leq z \leq 1}|(x_s)(z)|, \sup_{0 \leq z \leq 1}|(\tilde{x}_s)(z)|, |w(s)|\right)$$

$$L_a(R) := \sup\left\{\frac{|a(p,z) - a(\tilde{p},z)|}{|p - \tilde{p}|} : p \neq \tilde{p}, |p| \leq R, |\tilde{p}| \leq R, z \in [0,1]\right\}$$

Using (1.3), (3.1), (3.3), (3.21) and (3.22) we obtain:

$$|p(s) - \tilde{p}(s)| \leq 2L(R) \sup_{0 \leq z \leq 1}|x(s,z) - \tilde{x}(s,z)|, \text{ for almost all } s \in \left[t_0, \frac{t_1+\delta}{2}\right] \quad (3.27)$$

$$|v(t - c^{-1}z) - \tilde{v}(t - c^{-1}z)| \leq L(R)(t - c^{-1}z - t_0) \sup_{0 \leq s \leq t - c^{-1}z - t_0}|x(t - c^{-1}z, s) - \tilde{x}(t - c^{-1}z, s)|$$
$$+ L(R) \sup_{t - c^{-1}z - t_0 < s \leq 1}|x(t - c^{-1}z, s) - \tilde{x}(t - c^{-1}z, s)|,$$

for almost all $t \in \left(t_0, \min\left(\frac{t_1+\delta}{2}, t_0 + c^{-1}\right)\right)$ and $0 \leq z < c(t - t_0)$ (3.28)

We define for all $h \in \left(0, \min\left(\frac{\delta - t_1}{2}, c^{-1}\right)\right)$:

$$\Delta(h) := \sup_{t_0 < t \leq t_0 + h}\left(\sup_{c(t-t_0) < z \leq 1}|x(t,z) - \tilde{x}(t,z)|\right), \quad E(h) := \sup_{t_0 < t \leq t_0 + h}\left(\sup_{0 \leq z \leq c(t-t_0)}|x(t,z) - \tilde{x}(t,z)|\right) \quad (3.29)$$

Using (3.26), (3.27), (3.28) and definitions (3.29) we obtain for all $h \in \left(0, \min\left(\frac{\delta - t_1}{2}, c^{-1}\right)\right)$:

$$\Delta(h) \leq h\Gamma_1(R)\max(\Delta(h), E(h)) \quad (3.30)$$

$$E(h) \leq h\Gamma_1(R)\max(\Delta(h), E(h)) + h\Gamma_2(R)E(h) + \Gamma_2(R)\Delta(h) \quad (3.31)$$

where $\Gamma_1(R) := \exp(3M(R)c^{-1})(C + CRL_a(R)c^{-1} + 2RL_a(R))2L(R)$, $\Gamma_2(R) := \exp(M(R)c^{-1})L(R)$. Using (3.31) and assuming that $h \leq \frac{1}{2(\Gamma_1(R) + \Gamma_2(R))}$ we get $E(h) \leq \max(1, 2\Gamma_2(R))\Delta(h)$. In this case, inequality (3.30) implies $\Delta(h) \leq h\Gamma_1(R)\max(1, 2\Gamma_2(R))\Delta(h)$ which gives:

$$\Delta(h) = E(h) = 0, \text{ for all } h \in \left(0, \frac{1}{2}\min\left(c^{-1}, \frac{\delta - t_1}{2}, \frac{1}{2(\Gamma_1(R) + \Gamma_2(R))(1 + 2\Gamma_2(R))}\right)\right] \quad (3.32)$$

Implication (3.32), definitions (3.29) and the fact that $\sup_{t \in [0,t_0]}\sup_{z \in [0,1]}|(x_t)(z) - (\tilde{x}_t)(z)| = 0$ imply that $\sup_{t \in [0,t_0+h]}\sup_{z \in [0,1]}|(x_t)(z) - (\tilde{x}_t)(z)| = 0$. Since $t_0 \in (0, t_1]$ is arbitrarily close to $t_1$, we get $\sup_{t \in [0,t_1+h)}\sup_{z \in [0,1]}|(x_t)(z) - (\tilde{x}_t)(z)| = 0$.

Thus, the locally bounded mapping $[0, t_{\max}) \ni t \to x_t \in L^\infty([0,1]; \Re)$, which satisfies (1.3), (1.4), (3.3) for $t \in [0, t_{\max})$ a.e. and (3.9) for all $t \in [0, t_{\max})$ and almost all $z \in [0,1]$, which was found in Step 1 is unique.

<u>Step 3:</u> In this step we show that if $t_{\max} < +\infty$ then $\limsup_{t \to t_{\max}^-} \sup_{0 \leq z \leq 1} |(x_t)(z)| = +\infty$ and that there exist non-decreasing functions $P: \Re_+ \to \Re_+$, $Q: \Re_+ \to \Re_+$ such that for every $x_0 \in L^\infty([0,1];\Re)$, $y_0 \in L^\infty([0,1];\Re)$, $w \in L^\infty_{loc}([0,+\infty);W)$ estimate (3.10) holds.

Indeed, if $t_{\max} < +\infty$ for the unique solution $(p,v) \in L^\infty_{loc}([-c^{-1},t_{\max});\Re^{N+1})$ of (3.7) corresponding to input $w \in L^\infty_{loc}([0,+\infty);W)$ with initial condition $(p,v) \in L^\infty([-c^{-1},0);\Re^{N+1})$ satisfying $p(s) \equiv 0$ and (3.19), then Theorem 2.1 guarantees that $\limsup_{t \to t_{\max}^-} \|(p_t,v_t)\| = +\infty$. By virtue of inequality (3.2) and equations (1.3), (3.3), $\limsup_{t \to t_{\max}^-} \|(p_t,v_t)\| = +\infty$ cannot happen unless $\limsup_{t \to t_{\max}^-} \sup_{0 \leq z \leq 1} |(x_t)(z)| = +\infty$. Theorem 2.1 guarantees that there exist non-decreasing functions $\widetilde{P}: \Re_+ \to \Re_+$, $\widetilde{G}: \Re_+ \to \Re_+$ such that for every $v_0 \in L^\infty([-c^{-1},0);\Re)$, $\widetilde{v}_0 \in L^\infty([-c^{-1},0);\Re)$, $w \in L^\infty_{loc}([0,+\infty);W)$ it holds that:

$$\|p_t - \widetilde{p}_t\| + \|v_t - \widetilde{v}_t\| \leq \widetilde{G}(s(t)) \exp(\widetilde{P}(s(t))t) \|v_0 - \widetilde{v}_0\|, \text{ for all } t \in [0,\delta) \tag{3.33}$$

where $\delta = \min(t_{\max}(v_0,w), t_{\max}(\widetilde{v}_0,w))$, $s(t) := \max\left(\sup_{0 \leq \tau \leq t} \|(p_\tau,v_\tau)\|, \sup_{0 \leq \tau \leq t} \|(\widetilde{p}_\tau,\widetilde{v}_\tau)\|, \sup_{0 \leq \tau \leq t} \|w(\tau)\|\right)$ and $(p,v) \in L^\infty_{loc}([-c^{-1},t_{\max}(v_0,w));\Re^{N+1})$, $(\widetilde{p},\widetilde{v}) \in L^\infty_{loc}([-c^{-1},t_{\max}(\widetilde{v}_0,w));\Re^{N+1})$ are the unique solutions of (3.7) corresponding to input $w \in L^\infty_{loc}([0,+\infty);W)$ with initial conditions $(0,v_0) \in L^\infty([-c^{-1},0);\Re^{N+1})$ and $(0,\widetilde{v}_0) \in L^\infty([-c^{-1},0);\Re^{N+1})$, respectively.

Let $x_0 \in L^\infty([0,1];\Re)$, $\widetilde{x}_0 \in L^\infty([0,1];\Re)$, $w \in L^\infty_{loc}([0,+\infty);W)$ and consider the locally bounded mappings $[0,t_{\max}(x_0,w)) \ni t \to x_t \in L^\infty([0,1];\Re)$, $[0,t_{\max}(\widetilde{x}_0,w)) \ni t \to \widetilde{x}_t \in L^\infty([0,1];\Re)$, which satisfy (1.3), (1.4), (3.3) for $t \in [0,t_{\max})$ a.e. and (3.9) for all $t \in [0,t_{\max})$ and almost all $z \in [0,1]$ as well as the analogous equations with $\widetilde{x}$ replacing $x$ and $\widetilde{x}_0$ replacing $x_0$. The mappings $[0,t_{\max}(x_0,w)) \ni t \to x_t \in L^\infty([0,1];\Re)$, $[0,t_{\max}(\widetilde{x}_0,w)) \ni t \to \widetilde{x}_t \in L^\infty([0,1];\Re)$ are found by obtaining the solutions $(p,v) \in L^\infty_{loc}([-c^{-1},t_{\max}(v_0,w));\Re^{N+1})$, $(\widetilde{p},\widetilde{v}) \in L^\infty_{loc}([-c^{-1},t_{\max}(\widetilde{v}_0,w));\Re^{N+1})$ of (3.7) with input $w \in L^\infty_{loc}([0,+\infty);W)$ and initial conditions $(0,v_0) \in L^\infty([-c^{-1},0);\Re^{N+1})$, $(0,\widetilde{v}_0) \in L^\infty([-c^{-1},0);\Re^{N+1})$, respectively, and by using (3.5), (3.19) and

$$\widetilde{v}(t) = \exp\left(-\int_t^0 a(0,c(s-t))ds\right) \widetilde{x}_0(-ct), \text{ for almost all } t \in [-c^{-1},0) \tag{3.34}$$

Inequality (3.10) for appropriate non-decreasing functions $P: \Re_+ \to \Re_+$, $Q: \Re_+ \to \Re_+$ is a direct consequence of (3.26) with $t_0 = 0$, (3.33), (3.19), (3.34), (1.3), (3.3) and (3.2).

<u>Step 4:</u> In order to show that the mapping $t \to x_t \in L^\mu([0,1];\Re)$ is continuous on $[0,t_{\max})$ for every $\mu \in [1,+\infty)$, it suffices to use the classical semigroup property (identity (3.24)) and to show that for every $\mu \in [1,+\infty)$, $\varepsilon > 0$ there exists $\widetilde{\delta} > 0$ such that

$$\int_0^1 |x(t,z) - x_0(z)|^\mu dz \leq \varepsilon^\mu, \text{ for all } t \in [0,\widetilde{\delta}] \tag{3.35}$$

Since $C^0([0,1];\Re)$ is dense in $L^\mu([0,1];\Re)$ for every $\mu \in [1,+\infty)$ for every $\varepsilon > 0$ and $x_0 \in L^\infty([0,1];\Re)$ there exists $\tilde{x}_0 \in C^0([0,1];\Re)$ with $\int_0^1 |\tilde{x}_0(z) - x_0(z)|^\mu dz \leq (\varepsilon/3)^\mu$ and $\sup_{0 \leq z \leq 1} |\tilde{x}_0(z)| \leq \sup_{0 \leq z \leq 1} |x_0(z)|$. Moreover, there exists $\delta \in (0,1)$ such that $\sup_{0 \leq s \leq \delta} \sup_{s \leq z \leq 1} |\tilde{x}_0(z-s) - \tilde{x}_0(z)| \leq \frac{\varepsilon}{3}$. Using the inequality

$$\left( \int_\delta^1 |x_0(z-s) - x_0(z)|^\mu dz \right)^{1/\mu} \leq$$

$$\left( \int_s^1 |x_0(z-s) - \tilde{x}_0(z-s)|^\mu dz \right)^{1/\mu} + \left( \int_s^1 |\tilde{x}_0(z-s) - \tilde{x}_0(z)|^\mu dz \right)^{1/\mu} + \left( \int_s^1 |\tilde{x}_0(z) - x_0(z)|^\mu dz \right)^{1/\mu}$$

we obtain $\int_\delta^1 |x_0(z-s) - x_0(z)|^\mu dz \leq \varepsilon^\mu$ for all $s \in [0, \delta]$. Consequently, for every $\mu \in [1, +\infty)$, for every $\varepsilon > 0$ and $x_0 \in L^\infty([0,1];\Re)$ there exists $\delta \in (0,1)$ such that $\sup_{0 \leq s \leq \delta} \int_s^1 |x_0(z-s) - x_0(z)|^\mu dz \leq \varepsilon^\mu$.

Using the identity $\int_0^1 |x(t,z) - x_0(z)|^\mu dz = \int_0^{ct} |x(t,z) - x_0(z)|^\mu dz + \int_{ct}^1 |x(t,z) - x_0(z)|^\mu dz$ for all $t \in [0, c^{-1}]$, we get

$$\int_0^1 |x(t,z) - x_0(z)|^\mu dz \leq 2^\mu c R^\mu t + \int_{ct}^1 |x(t,z) - x_0(z)|^\mu dz, \text{ for all } t \in [0, \min(c^{-1}/2, t_{\max}/2)] \tag{3.36}$$

where $R = \sup_{0 \leq s \leq t_{\max}/2} \max\left( |p(s)|, |v(s)|, \sup_{0 \leq z \leq 1} |(x_s)(z)|, |w(s)| \right)$. Moreover, using (3.9) and definitions (3.13) we get for all $t \in [0, \min(c^{-1}/2, t_{\max}/2)]$ and almost all $z \in (ct, 1]$:

$$|x(t,z) - x_0(z)| \leq \exp(\tilde{M}(R)t)(C + \tilde{M}(R))Rt + |x_0(z-ct) - x_0(z)| \tag{3.37}$$

Combining (3.36) and (3.37) we obtain for all $t \in [0, \min(c^{-1}/2, t_{\max}/2)]$:

$$\int_0^1 |x(t,z) - x_0(z)|^\mu dz \leq 2^\mu c R^\mu t + 2^\mu \exp(\mu \tilde{M}(R)t)(C + \tilde{M}(R))^\mu R^\mu t^\mu + 2^\mu \int_{ct}^1 |x_0(z-ct) - x_0(z)|^\mu dz \tag{3.38}$$

Let $\delta \in (0,1)$ be such that $2^\mu \sup_{0 \leq s \leq \delta} \int_s^1 |x_0(z-s) - x_0(z)|^\mu dz \leq \frac{1}{3} \varepsilon^\mu$. It follows from the previous inequality in conjunction with (3.38) that (3.35) holds with $\tilde{\delta} = \min\left( \frac{1}{2c}, \frac{t_{\max}}{2}, \frac{\delta}{c}, \frac{\varepsilon}{3(2R)^\mu c}, \frac{\exp(-\tilde{M}(R)c^{-1})}{2R(C + \tilde{M}(R))} \left(\frac{\varepsilon}{3}\right)^{1/\mu} \right)$.

The proof is complete. ◁

Having proved Theorem 3.1, we are justified for arbitrary $x_0 \in L^\infty([0,1];\Re)$, $w \in L^\infty_{loc}([0,+\infty);W)$ to name the unique mapping $[0,t_{max}) \ni t \to x_t \in L^\infty([0,1];\Re)$ satisfying (1.3), (1.4), (3.3) for $t \in [0,t_{max})$ a.e. with $(p,v) \in L^\infty_{loc}([0,t_{max});\Re^{N+1})$, and (3.9) for all $t \in [0,t_{max})$ and almost all $z \in [0,1]$ as "*the solution $x_t$ of (1.2), (1.3), (1.4) with initial condition $x_0 \in L^\infty([0,1];\Re)$, corresponding to input $w \in L^\infty_{loc}([0,+\infty);W)$*", without explicit reference to the maximal existence time of the solution $t_{max} = t_{max}(x_0,w) \in (0,+\infty]$.

A number of system-theoretic properties can be proved by using the results of Theorem 3.1. More specifically, we study system (1.2), (1.3), (1.4) under assumptions (A1), (A2) and the following assumption:

**(A3)** $w = (d,u)$, *where $d \in D \subseteq \Re^{m_1}$ is a compact set, $u \in U \subseteq \Re^{m_2}$ is a locally compact set with $0 \in U$ and $m = m_1 + m_2$. Moreover, there exists $b \in K_\infty$ such that for every $x \in L^\infty([0,1];\Re)$, $d \in D$, $u \in U$, the following inequality holds:*

$$|G(d,u,x)| + \sum_{i=1}^{N} |K_i(d,u,x)| \leq b\left(\max\left(|u|, \sup_{0 \leq z \leq 1} |x(z)|\right)\right) \tag{3.39}$$

We are ready to give a list of properties for system (1.2), (1.3), (1.4) that are derived from the results of Theorem 3.1.

Property 1: Robustness of the equilibrium point.

Assumption (A3) guarantees that $0 \in L^\infty([0,1];\Re)$ is an equilibrium point for system (1.2), (1.3), (1.4), when $u \in U \subseteq \Re^{m_2}$ is zero for any $d \in D \subseteq \Re^{m_1}$. However, in order to study the robust stability properties of the equilibrium point of (1.2), (1.3), (1.4) we need a stronger assertion, namely that $0 \in L^\infty([0,1];\Re)$ is a robust equilibrium point for the system (1.2), (1.3), (1.4) with input $u \in U \subseteq \Re^{m_2}$ (see [10]).

**Theorem 3.2:** *Consider system (1.2), (1.3), (1.4) under assumptions (A1), (A2), (A3). The function $0 \in L^\infty([0,1];\Re)$ is a robust equilibrium point for the system (1.2), (1.3), (1.4) with input $u \in U \subseteq \Re^{m_2}$, i.e., for all $\varepsilon > 0$ and $T \geq 0$ there exists $\delta := \delta(\varepsilon,T) > 0$ such that for every $x_0 \in L^\infty([0,1];\Re)$, $d \in L^\infty([0,+\infty);D)$, $u \in L^\infty_{loc}([0,+\infty);U)$ with $\sup_{0 \leq z \leq 1}|x_0(z)| + \sup_{t \geq 0}|u(t)| < \delta$ the unique solution $t \to x_t \in L^\infty([0,1];\Re)$ of (1.2), (1.3), (1.4) with initial condition $x_0 \in L^\infty([0,1];\Re)$, corresponding to inputs $d \in L^\infty([0,+\infty);D)$, $u \in L^\infty_{loc}([0,+\infty);U)$ satisfies $\sup_{0 \leq z \leq 1}|(x_t)(z)| < \varepsilon$ for all $t \in [0,T]$.*

**Proof:** The proof is a direct consequence of Theorem 2.2: under assumption (A3) system (3.7) satisfies assumption (H3). Details are left to the reader.   ◁

Property 2: Robust Global Asymptotic Stability and Input-to-State Stability

Using the results contained in [10], we are in a position to define the notion of Robust Global Asymptotic Stability and Input-to-State Stability for system (1.2), (1.3), (1.4) under assumptions (A1), (A2), (A3).

**Definition 3.3:** *Consider system (1.2), (1.3), (1.4) under assumptions (A1), (A2), (A3) and assume that for every $x_0 \in L^\infty([0,1];\Re)$, $d \in L^\infty([0,+\infty);D)$, $u \in L^\infty_{loc}([0,+\infty);U)$, the solution $t \to x_t \in L^\infty([0,1];\Re)$ of (1.2), (1.3), (1.4) with initial condition $x_0 \in L^\infty([0,1];\Re)$, corresponding to inputs $d \in L^\infty([0,+\infty);D)$, $u \in L^\infty_{loc}([0,+\infty);U)$ exists for all $t \geq 0$. We say that system (1.2), (1.3), (1.4) is Input-to-State Stable (ISS) from the input $u \in U \subseteq \Re^{m_2}$ uniformly in $d \in D \subseteq \Re^{m_1}$, if there exists a continuous non-decreasing function $\gamma: \Re_+ \to \Re_+$ such that the following properties hold:*

**Robust Lagrange Stability:** *For every $\varepsilon > 0$ it holds that*

$$\sup\left\{\sup_{0\leq z\leq 1}|(x_t)(z)| - \sup_{0\leq s\leq t}\gamma(|u(s)|): t \geq 0, \sup_{0\leq z\leq 1}|x_0(z)| \leq \varepsilon, d \in L^\infty([0,+\infty);D), u \in L^\infty_{loc}([0,+\infty);U)\right\} < +\infty.$$

**Robust Lyapunov Stability:** *For every $\varepsilon > 0$ there exists $\delta := \delta(\varepsilon) > 0$ such that*

$$\sup\left\{\sup_{0\leq z\leq 1}|(x_t)(z)| - \sup_{0\leq s\leq t}\gamma(|u(s)|): t \geq 0, \sup_{0\leq z\leq 1}|x_0(z)| \leq \delta, d \in L^\infty([0,+\infty);D), u \in L^\infty_{loc}([0,+\infty);U)\right\} < \varepsilon.$$

**Uniform Robust Attractivity:** *For every $\varepsilon > 0$ and $R \geq 0$ there exists $\tau := \tau(\varepsilon, R) > 0$ such that*

$$\sup\left\{\sup_{0\leq z\leq 1}|(x_t)(z)| - \sup_{0\leq s\leq t}\gamma(|u(s)|): t \geq \tau, \sup_{0\leq z\leq 1}|x_0(z)| \leq R, d \in L^\infty([0,+\infty);D), u \in L^\infty_{loc}([0,+\infty);U)\right\} < \varepsilon.$$

*If $U = \{0\}$ then we say that $0 \in L^\infty([0,1];\Re)$ is (Uniformly) Robustly Globally Asymptotically Stable (RGAS) for (1.2), (1.3), (1.4).*

Lemma 2.1, Lemma 4.2, Theorem 2.2 in [10] and Theorems 3.1, 3.2 give us the following result.

**Theorem 3.4:** *Consider system (1.2), (1.3), (1.4) under assumptions (A1), (A2), (A3). Then the following statements are equivalent:*

**(a)** *System (1.2), (1.3), (1.4) is ISS from the input $u \in U \subseteq \Re^{m_2}$ uniformly in $d \in D \subseteq \Re^{m_1}$.*

**(b)** *There exists $\sigma \in KL$ and a continuous, non-decreasing function $\gamma: \Re_+ \to \Re_+$ such that for every $x_0 \in L^\infty([0,1];\Re)$, $d \in L^\infty([0,+\infty);D)$, $u \in L^\infty_{loc}([0,+\infty);U)$, the solution $t \to x_t \in L^\infty([0,1];\Re)$ of (1.2), (1.3), (1.4) with initial condition $x_0 \in L^\infty([0,1];\Re)$, corresponding to inputs $d \in L^\infty([0,+\infty);D)$, $u \in L^\infty_{loc}([0,+\infty);U)$ satisfies $\sup_{0\leq z\leq 1}|(x_t)(z)| \leq \sigma\left(\sup_{0\leq z\leq 1}|x_0(z)|, t\right) + \sup_{0\leq s\leq t}\gamma(|u(s)|)$ for all $t \geq 0$.*

**(c)** *System (3.7) is ISS from the input $u \in U \subseteq \Re^{m_2}$ uniformly in $d \in D \subseteq \Re^{m_1}$.*

*Furthermore, if $U = \{0\}$ then the following statement is equivalent to statements (a), (b):*

**(d)** *System (1.2), (1.3), (1.4) is Robustly Forward Complete, i.e., for every $T \geq 0$ and $R \geq 0$ it holds that $\sup\left\{\sup_{0\leq z\leq 1}|(x_t)(z)|: t \in [0,T], \sup_{0\leq z\leq 1}|x_0(z)| \leq R, d \in L^\infty([0,+\infty);D)\right\} < +\infty$, where $x_t \in L^\infty([0,1];\Re)$ denotes the solution of (1.2), (1.3), (1.4) with initial condition $x_0 \in L^\infty([0,1];\Re)$, corresponding to input $d \in L^\infty([0,+\infty);D)$, and the property of Uniform Robust Attractivity of Definition 3.3 holds.*

Property 3: Lyapunov Characterization of RGAS

Theorem 3.4 in [10] and the results of Theorems 3.1, 3.2, 3.4 allow us to obtain a complete Lyapunov characterization of the RGAS property for system (1.2), (1.3), (1.4).

**Theorem 3.5:** *Consider system (1.2), (1.3), (1.4) under assumptions (A1), (A2), (A3) and assume that $U = \{0\}$. The equilibrium point $0 \in L^\infty([0,1]; \Re)$ is RGAS for (1.2), (1.3), (1.4) if and only if there exists a functional $V: L^\infty([0,1]; \Re) \to \Re_+$, a non-decreasing function $Q: \Re_+ \to \Re_+$ and functions $a_1, a_2 \in K_\infty$ such that the following inequalities hold:*

$$a_1\left(\sup_{0 \leq z \leq 1} |x(z)|\right) \leq V(x) \leq a_2\left(\sup_{0 \leq z \leq 1} |x(z)|\right), \text{ for all } x \in L^\infty([0,1]; \Re) \quad (3.40)$$

$$|V(x) - V(y)| \leq Q\left(\max\left(\sup_{0 \leq z \leq 1} |x(z)|, \sup_{0 \leq z \leq 1} |y(z)|\right)\right) \sup_{0 \leq z \leq 1} |x(z) - y(z)|, \text{ for all } x, y \in L^\infty([0,1]; \Re) \quad (3.41)$$

*Moreover, for every $x_0 \in L^\infty([0,1]; \Re)$, $d \in L^\infty([0,+\infty); D)$, the solution $x_t \in L^\infty([0,1]; \Re)$ of (1.2), (1.3), (1.4) with initial condition $x_0 \in L^\infty([0,1]; \Re)$, corresponding to input $d \in L^\infty([0,+\infty); D)$ satisfies:*

$$V(x_t) \leq \exp(-t) V(x_0), \text{ for all } t \geq 0 \quad (3.42)$$

Again, it should be pointed out that inequality (3.41) guarantees that the functional $V: L^\infty([0,1]; \Re) \to \Re_+$ is Lipschitz on bounded sets of the state space $L^\infty([0,1]; \Re)$. However, inequality (3.41) does not guarantee Frechet differentiability of the functional $V: L^\infty([0,1]; \Re) \to \Re_+$ nor that the limit $\lim_{t \to 0^+} t^{-1}(V(x_t) - V(x_0))$ exists for the solution $x_t$ of (1.2), (1.3), (1.4) with initial condition $x_0 \in L^\infty([0,1]; \Re)$, corresponding to input $d \in L^\infty([0,+\infty); D)$. Notice that inequality (3.42) guarantees that for every $x_0 \in L^\infty([0,1]; \Re)$, $d \in L^\infty([0,+\infty); D)$, the solution $x_t$ of (1.2), (1.3), (1.4) with initial condition $x_0 \in L^\infty([0,1]; \Re)$, corresponding to input $d \in L^\infty([0,+\infty); D)$ satisfies $\limsup_{t \to 0^+} t^{-1}(V(x_t) - V(x_0)) \leq -V(x_0)$.

Property 4: Sufficient Conditions for Stability Properties

Theorem 3.5 is not the most convenient way of proving RGAS or ISS for (1.2), (1.3), (1.4). For practical purposes we can use the following result, which is a direct consequence of Theorem 2.6 and the equivalent description of system (1.2), (1.3), (1.4) by means of system (3.7). Its proof is omitted.

**Theorem 3.6:** *Consider system (1.2), (1.3), (1.4) under assumptions (H1), (H2), (H3). Assume that there exists a continuous, positive definite and radially unbounded function $W: \Re^{N+1} \to \Re_+$, a continuous, non-decreasing function $\gamma: \Re_+ \to \Re_+$ and a constant $\lambda \in (0,1)$ such that the following inequality holds for all $(p,v) \in L^\infty([-c^{-1}, 0); \Re^{N+1})$, $d \in D$, $u \in U$:*

$$W(F(d, u, p, v)) \leq \lambda \sup_{-c^{-1} \leq s < 0} W(p(s), v(s)) + \gamma(|u|) \quad (3.43)$$

where $F(d,u,p,v) := (F_1(d,u,p,v),...,F_{N+1}(d,u,p,v))$ and $F_i : D \times U \times L^\infty([-c^{-1},0];\Re^{N+1}) \to \Re$ ($i=1,...,N+1$) are defined by (3.8). Then system (1.2), (1.3), (1.4) is ISS from the input $u \in U \subseteq \Re^{m_2}$ uniformly in $d \in D \subseteq \Re^{m_1}$.

The use of Theorem 3.6 is illustrated by the following example.

**Example 3.7:** Consider system (1.2), (1.3), (1.4), where $c=1$, $a(p,z) \equiv 0$, $N=1$, $U=\{0\}$, $K_1(d,u,x) = d\int_0^1 x(z)dz$, $G(d,u,x) \equiv 0$, $g_1(z) = g$, $g \in \Re$ is a constant and $D := [-1,1] \subset \Re$, namely the system

$$\frac{\partial x}{\partial t}(t,z) + \frac{\partial x}{\partial z}(t,z) = gd(t)\int_0^1 x(t,z)dz, \quad x(t,0) = 0, \quad |d(t)| \leq 1 \tag{3.44}$$

For this example, system (3.7) is given by:

$$p(t) = d(t)\int_{t-1}^t v(s)ds + gd(t)\int_{t-1}^t (1+s-t)p(s)ds \tag{3.45}$$
$$v(t) = 0$$

and the solution of (3.44) is related to the solution of (3.45) by means of the equations:

$$x(t,z) = v(t-z) + g\int_{t-z}^t p(w)dw, \quad p(t) = \int_0^1 x(t,z)dz \tag{3.46}$$

Clearly, assumptions (H1), (H2), (H3) are satisfied for system (3.45) and assumptions (A1), (A2), (A3) are satisfied for system (3.44). In this example we show that the equilibrium point $0 \in L^\infty([0,1];\Re)$ is RGAS for system (3.44) provided that:

$$|g| < 2 \tag{3.47}$$

Indeed, we apply Theorem 3.6 with $F(d,u,p,v) := \left(d\int_{-1}^0 v(s)ds + gd\int_{-1}^0 (1+s)p(s)ds, 0\right)'$ and $W(p,v) = |p| + k|v|$ where $k > \frac{2}{2-|g|}$ is a constant. Notice that:

$$W(F(d,u,p,v)) = |d|\left|\int_{-1}^0 v(s)ds + g\int_{-1}^0 (1+s)p(s)ds\right| \leq \left|\int_{-1}^0 v(s)ds + g\int_{-1}^0 (1+s)p(s)ds\right|$$
$$\leq \sup_{-1 \leq s < 0} |v(s)| + \frac{1}{2}|g| \sup_{-1 \leq s < 0} |p(s)| \leq \left(\frac{1}{2}|g| + k^{-1}\right) \sup_{-1 \leq s < 0} (|p(s)| + k|v(s)|) = \left(\frac{1}{2}|g| + k^{-1}\right) \sup_{-1 \leq s < 0} W(p(s), v(s))$$

Consequently, (3.43) holds with $\lambda := \frac{1}{2}|g| + k^{-1} < 1$ (with $\gamma \equiv 0$). The condition (3.47) is sharp: for $g=2$ and $d(t) \equiv 1$, system (3.44) admits the family of constant solutions $x(t,z) = Qz$ for $z \in [0,1]$, where $Q \in \Re$ and consequently the equilibrium point $0 \in L^\infty([0,1];\Re)$ is not RGAS. ◁

# 4. An Illustrative Example

Consider system (1.2), (1.3), (1.4), where $c=1$, $a(p,z) \equiv 0$, $N=1$, $U = \Re$, $K_1(u,x) = x(1)$, $g_1(z) = g$, $g \in \Re$ is a constant and $d \in D$ is irrelevant. This problem was studied in [12] and a boundary feedback control of the form:

$$G(u,x) = \int_0^1 k(z)x(z)dz \qquad (4.1)$$

was designed, where the kernel $k \in C^1([0,1];\Re)$ was explicitly given and guaranteed finite time stability for the corresponding closed-loop system. Here we will provide an alternative methodology for the design of a robust stabilizer for system (1.2), (1.3), (1.4) with $G(u,x) = u \in \Re$, namely, the system

$$\frac{\partial x}{\partial t}(t,z) + \frac{\partial x}{\partial z}(t,z) = g\,x(t,1), \quad x(t,0) = u(t) \qquad (4.2)$$

For this example, system (3.7) is given by:

$$p(t) = v(t-1) + g\int_{-1}^{0} p(t+s)ds$$
$$v(t) = u(t) \qquad (4.3)$$

and the solution of (4.2) is related to the solution of (4.3) by means of the equations:

$$x(t,z) = v(t-z) + g\int_{t-z}^{t} p(w)dw, \quad p(t) = x(t,1) \qquad (4.4)$$

Clearly, assumptions (H1), (H2), (H3) are satisfied for system (4.3) and assumptions (A1), (A2), (A3) are satisfied for system (4.2). Equations (4.3) imply that the following differential equation holds for almost all $t \geq 0$ for which the solution of (4.3) exists:

$$\frac{d}{dt}(p(t) - v(t-1)) = g(p(t) - v(t-1)) - gp(t-1) + gv(t-1) \qquad (4.5)$$

Consequently, the following equation holds for all $t \geq t_0 \geq 0$ for which the solution of (4.3) exists:

$$p(t) = v(t-1) + \exp(g(t-t_0))(p(t_0) - v(t_0-1))$$
$$- g\int_{t_0}^{t} \exp(g(t-s))p(s-1)ds + g\int_{t_0}^{t} \exp(g(t-s))v(s-1)ds \qquad (4.6)$$

Using (4.6) with $t_0 = 0$, we can show that for every initial condition and for every $u \in L^{\infty}_{loc}(\Re_+;\Re)$ the solution of (4.3) is bounded for all $t \in [0,1)$: therefore, Theorem 2.2 implies that the solution exists for all $t \in [0,1)$. Using induction we conclude that for every initial condition and for every $u \in L^{\infty}_{loc}(\Re_+;\Re)$ the corresponding solution of (4.3) exists for all $t \geq 0$. Moreover, using (4.3) and (4.6) with $t_0 = t-1$ we conclude that system (4.3) satisfies the following equations for $t \geq 1$ a.e.:

$$p(t) = v(t-1) + g \int_{t-2}^{t-1} (\exp(g) - \exp(g(t-s-1)))p(s)ds + g \int_{t-2}^{t-1} \exp(g(t-s-1))v(s)ds \quad (4.7)$$

$$v(t) = u(t)$$

It is clear that the feedback law:

$$u(t) = -g \int_{t-1}^{t} (\exp(g) - \exp(g(t-s)))p(s)ds - g \int_{t-1}^{t} \exp(g(t-s))v(s)ds \quad (4.8)$$

guarantees that the closed-loop system (4.3) with (4.8) satisfies for $t \geq 1$ a.e.:

$$p(t) = 0$$
$$v(t) = -g \int_{-1}^{0} (\exp(-gw)v(t+w) + (\exp(g) - \exp(-gw))p(t+w))dw \quad (4.9)$$

Using (4.3) (and the equation $v(t) = p(t+1) - g \int_{-1}^{0} p(t+1+s)ds$) we can conclude from (4.9) that $v(t) = 0$ for $t \geq 1$.

At this point, we should emphasize that we have achieved finite-time stabilization of the system (4.3) by means of the feedback law (4.8). However, we have not achieved stabilization of system (4.2) by means of a feedback law of the form (4.1). This remains to be shown. Indeed, using (4.4) we get from (4.8):

$$u(t) = -g \int_{t-1}^{t} (\exp(g) - \exp(g(t-s)))p(s)ds - g \int_{t-1}^{t} \exp(g(t-s))v(s)ds$$

$$= -g \int_{t-1}^{t} (\exp(g) - \exp(g(t-s)))p(s)ds - g \int_{0}^{1} \exp(gz)v(t-z)dz$$

$$= -g \int_{t-1}^{t} (\exp(g) - \exp(g(t-s)))p(s)ds + g^2 \int_{0}^{1} \exp(gz) \int_{t-z}^{t} p(s)dsdz - g \int_{0}^{1} \exp(gz)x(t,z)dz$$

Integrating by parts we get:

$$g^2 \int_{0}^{1} \exp(gz) \int_{t-z}^{t} p(s)dsdz = g \int_{0}^{1} \frac{d}{dz}(\exp(gz)) \int_{t-z}^{t} p(s)dsdz$$

$$= g \int_{0}^{1} \frac{d}{dz}\left(\exp(gz) \int_{t-z}^{t} p(s)ds\right)dz - g \int_{0}^{1} \exp(gz) \frac{d}{dz}\left(\int_{t-z}^{t} p(s)ds\right)dz$$

$$= g \exp(g) \int_{t-1}^{t} p(s)ds - g \int_{0}^{1} \exp(gz)p(t-z)dz = g \int_{t-1}^{t} (\exp(g) - \exp(g(t-s)))p(s)ds$$

and consequently we obtain $u(t) = -g \int_{0}^{1} \exp(gz)x(t,z)dz$, which is a boundary feedback law of the form (4.1) with $k(z) = -g \exp(gz)$. Consequently, using (4.4) we conclude that the feedback law (4.1) with $k(z) = -g \exp(gz)$ achieves finite-time stabilization of system (4.2). Indeed, we can prove the existence of a constant $M \geq 1$ so that the following estimate holds for the solution of the closed-loop system (4.3) with (4.8):

$$\sup_{t-1\leq s<t}|v(s)| + \sup_{t-1\leq s<t}|p(s)| \leq M\sigma(t)\left(\sup_{-1\leq s<0}|v(s)| + \sup_{-1\leq s<0}|p(s)|\right), \text{ for all } t \geq 0 \qquad (4.10)$$

where $\sigma(t) = \begin{cases} 1 & \text{if } t \in [0,2) \\ 0 & \text{if } t \geq 2 \end{cases}$. Using (4.4) we get $\sup_{0<z\leq 1}|x(t,z)| \leq \sup_{t-1\leq s<t}|v(s)| + |g|\sup_{t-1\leq s<t}|p(s)|$. The previous inequality, estimate (4.10) and the fact that (4.4) must hold, allow us (by selecting $p(s) \equiv 0$ for $s \in [-1,0)$, which by virtue of (4.4) gives $v(-q) = x(0,q)$ for $q \in (0,1]$) to obtain the estimate:

$$\sup_{0<z\leq 1}|x(t,z)| \leq (1+|g|)M\sigma(t)\sup_{0<z\leq 1}|x(0,z)|, \text{ for all } t \geq 0 \qquad (4.11)$$

The analysis presented so far allowed us to obtain some different features from the analysis in [12]:

1) The implementation of the feedback law: It follows from (4.4) that the feedback law (4.8) can be implemented by using the following distributed delay feedback law:

$$u(t) = -g\int_{t-1}^{t}(\exp(g) - \exp(g(t-s)))x(s,1)ds - g\int_{t-1}^{t}\exp(g(t-s))x(s,0)ds \qquad (4.12)$$

The feedback law (4.12) shows that we do not need to measure the whole state profile: we only need to measure $x(t,z)$ at two different space points $z=0$ and $z=1$.

2) The issue of the state space: It should be emphasized that the feedback law $u(t) = -g\int_{0}^{1}\exp(gz)x(t,z)dz$ is exactly the same with that obtained in [12]. However, we have proved that estimate (4.11) holds for initial conditions $x(0,z) = x_0(z)$ for $z \in [0,1]$ with $x_0 \in L^{\infty}([0,1];\Re)$. Moreover, no compatibility condition (e.g., $\lim_{z\to 0^+} x_0(z) = -g\int_0^1 \exp(gz)x_0(z)dz$) is required to hold.

3) Control actuator errors: The implementation of the feedback law $u(t) = -g\int_0^1 \exp(gz)x(t,z)dz$ may result to the equation:

$$u(t) = -g\int_0^1 \exp(gz)x(t,z)dz + w(t), \text{ for } t \geq 0 \qquad (4.13)$$

where $w \in L^{\infty}(\Re_+;\Re)$ is the control actuator error. Using the framework described in the present work, we are in a position to study the closed-loop system (4.2) with (4.13). In this case, we get the integral delay system:

$$p(t) = v(t-1) + g\int_{-1}^{0} p(t+s)ds$$
$$v(t) = w(t) - g\int_{-1}^{0}(\exp(g) - \exp(-gs))p(t+s)ds - g\int_{-1}^{0}\exp(-gs)v(t+s)ds \qquad (4.14)$$

with $w \in L^{\infty}(\Re_+;\Re)$ as input. Using (4.7), we get for all $t \geq 1$:

$$p(t) = w(t-1) \quad , \quad v(t) = w(t) - g\int_{-1}^{0} w(t+s)ds \qquad (4.15)$$

which can give us the existence of constants $M \geq 1$, $\gamma \geq 0$ such that the following estimate holds:

$$\sup_{t-1 \leq s < t} |v(s)| + \sup_{t-1 \leq s < t} |p(s)| \leq M\sigma(t)\left(\sup_{-1 \leq s < 0} |v(s)| + \sup_{-1 \leq s < 0} |p(s)|\right) + \gamma \sup_{0 \leq s \leq t} |w(s)|, \text{ for all } t \geq 0 \quad (4.16)$$

where $\sigma(t) = \begin{cases} 1 & \text{if } t \in [0,2) \\ 0 & \text{if } t \geq 2 \end{cases}$. Estimate (4.16) shows that the ISS property holds for system (4.14). Estimate (4.16) shows (using exactly the same arguments as above) that the following estimate holds for the solution of (4.2) with (4.13):

$$\sup_{0 < z \leq 1} |x(t,z)| \leq (1+|g|)M\sigma(t)\sup_{0 < z \leq 1}|x(0,z)| + (1+|g|)\gamma \sup_{0 \leq s \leq t}|w(s)|, \text{ for all } t \geq 0 \quad (4.17)$$

which is exactly the ISS property for the closed-loop system (4.2) with (4.13).

Having presented the analysis of the example, we are in a position to point out certain shortcomings of the proposed methodology of the conversion of a single FOH-PDE to a system described by IDEs:

i) The major shortcoming of the above analysis is that the researcher studies a different system (an integral delay system) from the original one (a system described by first-order hyperbolic pdes). This is important, because a stabilizing feedback law for (4.3) may not be equivalent to a state feedback law for (4.2): it may involve delay terms. Therefore, a stabilizing feedback law for (4.3), when expressed in the original "coordinates" may not give us a closed-loop system of the form (4.2) but rather **a system described by a single first-order hyperbolic PDE with delays**.

ii) Another shortcoming of the above analysis is the specific form (1.2), (1.3), (1.4) for which the conversion to a system described by IDEs can be used. Not all systems described by first-order hyperbolic pdes can be expressed by (1.2), (1.3), (1.4): the requirement that the functionals $K_i : W \times L^\infty([0,1]; \Re) \to \Re$ ($i = 1,...,N$) are not allowed to depend on $z \in [0,1]$ is restrictive. There are systems which are studied in [12] but cannot be studied in the proposed framework. However, one must bear in mind that a preliminary integral transformation can be applied so that the transformed system is of the form (1.2), (1.3), (1.4).

## 5. Concluding Remarks

This paper establishes the equivalence between systems described by a single first-order hyperbolic partial differential equation and systems described by integral delay equations. System-theoretic results are provided for both classes of systems:
1) the Boundedness-Implies-Continuation property and the property of Lipschitz dependence on initial conditions (Theorem 2.1 and Theorem 3.1),
2) the robustness of an equilibrium point for systems with external inputs (Theorem 2.2 and Theorem 3.2),
3) characterizations of the Input-to-State Stability and Robust Global Asymptotic Stability (Theorem 2.4 and Theorem 3.4),
4) converse Lyapunov results (Theorem 2.5 and Theorem 3.5)
5) sufficient conditions for Input-to-State Stability and Robust Global Asymptotic Stability (Theorem 2.6 and Theorem 3.6).

The proposed framework can allow the study of discontinuous solutions for nonlinear systems described by a single first-order hyperbolic partial differential equation under the effect of measurable inputs acting on the boundary and/or on the differential equation. This aspect is important from a control-theoretic point of view because all systems are subject to disturbances (modelling errors, control actuator errors and measurement errors), which are typically modeled by measurable inputs. An illustrative example shows that the conversion of a system described by a single first-order hyperbolic partial differential equation to an integral delay system can simplify considerably the solution of the corresponding robust feedback stabilization problem.